\documentclass[12pt,a4paper]{amsart}

\usepackage{amsmath}
\usepackage{amssymb}
\usepackage{latexsym}
\usepackage{xypic}
\usepackage{amscd,amsthm}

\usepackage{textcomp}
\usepackage{graphicx}
\usepackage{psfrag}
\usepackage{epic}
\usepackage{eepic}
\usepackage[matrix,arrow,curve]{xy}
\input{epsf}

\def\C{{\mathbb C}}
\def\N{{\mathbb N}}

\def\P{{\mathbb P}}
\def\Q{{\mathbb Q}}
\def\R{{\mathbb R}}
\def\Z{{\mathbb Z}}
\def\A{{\mathbb A}}

\def\codim{{\rm codim}}

\def\deg{{\rm deg}}

\newtheorem*{theorem*}{Theorem}
\newtheorem*{corollary*}{Corollary}
\newtheorem{theorem}{Theorem}[section]
\newtheorem{definition}[theorem]{Definition}
\newtheorem{remark}[theorem]{Remark}

\newtheorem{corollary}[theorem]{Corollary}
\newtheorem{proposition}[theorem]{Proposition}

\begin{document}
\title[$SL(2)$-equivariant flips]{On the geometry of 
$SL(2)$-equivariant flips}

\author[Victor Batyrev]{Victor Batyrev}
\address{Mathematisches Institut, Universit\"at T\"ubingen, 
Auf der Morgenstelle 10, 72076 T\"ubingen, Germany}
\email{ victor.batyrev@uni-tuebingen.de}

\author[Fatima Haddad]{Fatima Haddad}
\address{Mathematisches Institut, Universit\"at T\"ubingen, 
Auf der Morgenstelle 10, 72076 T\"ubingen, Germany}
\email{fatima@everest.mathematik.uni-tuebingen.de}

\begin{abstract}
\vspace{-0.5ex}
In this paper, we show that  any  $3$-dimensional 
normal affine quasihomogeneous $SL(2)$-variety 
can be described as a categorical quotient of a $4$-dimensional 
affine hypersurface. Moreover, 
we show that the Cox ring of an arbitrary $3$-dimensional 
normal affine quasihomogeneous $SL(2)$-variety has a unique 
defining equation. This allows us to construct  $SL(2)$-equivariant
flips by different GIT-quotients of hypersurfaces.  Using 
the theory of spherical varieties, 
we describe  $SL(2)$-flips by means of $2$-dimensional 
colored cones.  
\end{abstract}

\dedicatory{Dedicated to Ernest Borisovich Vinberg on 
the occasion of his 70th birthday}

\maketitle

\section*{Introduction}

Let $X$, $X^-$ and $X^+$ be normal quasiprojective $3$-dimensional 
algebraic varieties over $\C$ \footnote{All results of our  paper 
are valid for  algebraic varieties defined over an arbitrary 
algebraically closed field $K$ of characteristic $0$, but 
for simplicity we consider 
only the case $K=\C$.}. A flip  is 
a diagram 
$$
\xymatrix{ X^- \ar@{-->}[rr] \ar[dr]_{\varphi^-}& &  X^+ \ar[dl]^{\varphi^+}\\
& X & }
$$ 
in which $X^-$ and $X^+$ are $\Q$-factorial and   
 $\varphi^-\, : \, X^- \to X$, 
 $\varphi^+\, : \, X^+ \to X$  are projective birational morphisms 
contracting  
finitely many rational curves to an isolated singular point $p \in X$. 
Moreover,  the anticanonical divisor $-K_{X^-}$ and the canonical 
divisor $K_{X^+}$ 
are relatively ample over $X$. In this case, 
the singulartity at $p$ is not $\Q$-factorial 
(even not $\Q$-Gorenstein).  
A good understanding of flips is important from 
the point of view of $3$-dimensional birational geometry (see e.g. 
\cite{CKM88}, or \cite{KM98}).  
We remark that if $X$ is affine, then 
the quasiprojective varieties $X^-$ and $X^+$ can be 
obtained from $X$ as follows: 
\[ X^- := {\rm Proj}\, \bigoplus_{n \geq 0}  
\Gamma(X,{\mathcal O}_X(-nK_X)), \;\; 
X^+ := {\rm Proj}\, \bigoplus_{n \geq 0} 
\Gamma(X, {\mathcal O}_X(nK_X)). \]   

Simplest examples of flips can be constructed from 
$3$-dimensional affine toric varieties $X_{\sigma}$
which are  categorical quotients of $\C^4$ modulo $\C^*$-actions by 
diagonal matrices $$diag(t^{-n_1},t^{-n_2}, t^{n_3}, 
t^{n_4}), \;\;  t\in \C^*$$ where $n_1,n_2,
n_3,n_4$ are 
positive integers satisfying the condition  
$n_1 + n_2 < n_3+n_4$ \cite{Da83,R83}. 
In this case, the quasiprojective 
toric varieties $X_{\sigma}^-$ and $X_{\sigma}^+$ 
correspond to two different simplicial 
subdivisions
of a $3$-dimensional cone $\sigma$  generated by 
$4$ lattice vectors $v_1, v_2, v_3, v_4$ satisfying the relation 
$n_1 v_1 + n_2 v_2 = n_3 v_3 + n_4 v_4$.       

Another point of view on flips comes from the Geometric Invariant 
Theory (GIT) which describes a flip diagram as 
$$
\xymatrix{ Y^{ss}(L^-)/\!\!/ G \ar@{-->}[rr] \ar[dr]_{\varphi^-}& &   
Y^{ss}(L^+)/\!\!/ G \ar[dl]^{\varphi^+}\\
& Y^{ss}(L^0)/\!\!/ G  & }
$$ 
for some three $G$-linearized ample line bundles 
$L^-$, $L^+$, $L^0$ on a $4$-dimensional variety $Y$ (see e.g. 
\cite{Th96}). Here 
\[ Y^{ss}(L) := \{ y \in Y\; : s(y) \neq 0 \;\; {\rm  for} \;{\rm  
some}\; s \in \Gamma(Y, L^{\otimes n})^G  \;{\rm and} \;\; {\rm  for} \;{\rm  
some}\; n >0  \} \]
denotes the subset 
of semistable  points in $Y$ 
with respect to the $G$-linearized ample line 
bundle $L$  and $Y^{ss}(L)/\!\!/G$ denotes 
the categorical quotient which can be identified with 
\[ {\rm Proj}\,  \bigoplus_{n \geq 0} \Gamma(Y, L^{\otimes n})^G. \]
In the above toric case,
we have $Y = \C^4$, $G \cong \C^*$ and   $L^-$, $L^+$, $L^0$
are different $G$-linearizations of the trivial line bundle 
over $\C^4$.   
A classification of $3$-dimensional  flips in case when 
$Y \subset \C^5$ is hypersurface and 
both varieties  $Y^{ss}(L^-)/\!\!/ G$ and 
$ Y^{ss}(L^+)/\!\!/ G$ have 
at worst terminal singularities was 
considered by Brown in \cite{Br99}. 
\medskip

The purpose of this paper is to investigate another class 
of quasihomogeneous varieties.  We  
give a geometric description of 
$SL(2)$-equivariant flips 
in the case when $X$ is an arbitrary  
singular  normal  affine quasihomogeneous 
$SL(2)$-variety.  
 
It follows form a recently result of Gaifullin \cite{Ga08} that a 
$3$-dimensional toric 
variety $X_{\sigma}$ associated with a $3$-dimensional 
cone $\sigma = \R_{\geq 0}v_1 + \cdots +  \R_{\geq 0}v_4$ 
is quasihomogeneous with respect to 
a $SL(2)$-action if and only if $n_1 = n_2$ and $n_3 = n_4$. 
One of  such  toric  $SL(2)$-equivariant flips $(n_1=n_2 =1, \; 
n_3 = n_4 = n >1$)  
was described in detail 
in \cite[Example 2.7]{KM98}. However, there exist many normal affine 
$3$-dimensional quasihomogeneous $SL(2)$-varieties which are 
not toric.  
\medskip

According to Popov \cite{P73}, every  normal  affine quasihomogeneous 
$SL(2)$-variety $E$ is uniquely determined by a pair of numbers  
$(h, m) \in \{\Q \cap (0,1] \}\times \N$ (we denote this variety by 
$E_{h,m}$). 
Let $h = p/q \leq 1$ $(g.c.d.(p,q)=1)$. We define  $k:=g.c.d.(q-p,m)$,   
\[   a := \frac{m}{k},\;\;  b := \frac{(q-p)}{k}.\]
In Section 1 we show that    
the affine $SL(2)$-variety $E_{h,m}$ is isomorphic to the 
catego\-ri\-cal quotient of the 
hypersurface $H_b \subset \C^5$  defined by the equation 
\[ Y_0^b = X_1X_4 - X_2X_3, \] 
modulo the action of the diagonalizable group  
$G \cong \C^*\times \mu_a$, where 
$\C^*$ acts by $$diag(t^k,t^{-p},t^{-p},t^q,t^q),\;\;  t \in \C^*$$ and 
$\mu_a = \langle \zeta_a \rangle$, $\zeta_a = e^{2\pi i/a}$ 
acts by 
 $diag(1,\zeta_a^{-1}, \zeta_a^{-1}, \zeta_a, \zeta_a )$. 
Here we consider the $SL(2)$-action on $H_b$  
 induced by the trivial action on the coordinate $Y_0$ and by left 
multiplication on the coordinates $X_1, X_2, X_3, X_4$: 
 \[  \left( \begin{pmatrix} X & Y \\ Z & W \end{pmatrix}, 
\begin{pmatrix} X_1 & X_3 \\ X_2 & X_4 \end{pmatrix} \right) \mapsto 
\begin{pmatrix} X & Y \\ Z & W \end{pmatrix}\cdot 
 \begin{pmatrix} X_1 & X_3 \\ X_2 & X_4 \end{pmatrix}, \; \;\; 
 \;\; \begin{pmatrix} X & Y \\ Z & W \end{pmatrix} \in 
SL(2). \]
This $SL(2)$-action commutes with the $G$-action and descends to the 
categorial quotient $H_b /\!\!/G \cong E_{h,m}$.  
In this way,  
we obtain  a very simple description of the affine $SL(2)$-variety 
$E_{h,m}$ which seems to be overlooked in the literature.
 \medskip

In Section 2 we consider  the notion of the 
total coordinate ring (or Cox ring) 
of an  algebraic variety 
$X$ with  a finitely generated divisor class group ${\rm Cl}(X)$ 
(see e.g. \cite{Ar08,H08}. 
These rings naturally appear  in some questions related to 
Del Pezzo surfaces and homogeneous spaces of 
algebraic groups  (see \cite{BP04}). 
Using  results from  Section 1,  we show that  
the Cox ring  of $E_{h,m}$ 
is isomorphic to 
\[ \C[Y_0, X_1, X_2, X_3, X_4] / (  Y^b_0 - X_1 X_4 + X_2X_3). \]
Some similar examples of algebraic varieties whose Cox ring 
is defined by a unique equation were considered in \cite{BH07}. 
We remark that our  result provides an alternative 
proof of a  criterion  of Gaifullin \cite{Ga08}: 
$E_{h,m}$ is toric if and only 
if $b =1$, or  if only if  $(q-p)$ divides $m$. 
One can use our description  of the total coordinate ring 
of $E_{h,m}$  as a good illustration of  more general 
resent results of Brion on the total coordinate ring of spherical 
varieties \cite{B07}.

\medskip

In Section 3 we  describe the  quasiprojective varieties $E_{h,m}$
$E_{h,m}^+$,  $E_{h,m}^-$ in the $SL(2)$-equivariant flip diagram 
\[
\xymatrix{ E_{h,m}^- = H_b^{ss}(L^-)/\!\!/G  \ar@{-->}[rr] \ar[dr]_{\varphi^-}& &
 E_{h,m}^+ = H_b^{ss}(L^+)/\!\!/G  \ar[dl]^{\varphi^+}\\
&E_{h,m} = H_b^{ss}(L^0)/\!\!/G  & }
\]
by different GIT-quotients of the hypersurface $H_b$, where $L^0$, 
$L^-$, $L^+$  are linearizations of the  trivial line 
bundle over $H_b$ corresponding to the trivial character $\chi^0$ 
and some nontrivial characters $\chi^-$, $\chi^+$ 
of $G$. In fact, we can 
say something more about  $E_{h,m}^-$ and  $E_{h,m}^+$: 
there exist two affine normal toric 
surfaces $S^-  \subset E_{h,m}$  and $S^+ \subset  E_{h,m}$ 
which are closures in $E_{h,m}$ 
of two orbits of a Borel subgroup $B \subset SL(2)$ such that 
\[  E_{h,m}^- \cong SL(2) \times_B S^-, 
\;\; E_{h,m}^+ \cong SL(2) \times_B S^+.\]
In particular,  both varieties $E_{h,m}^+$ and   
$E_{h,m}^-$ have at worst log-terminal 
toroidal singularities. 
The birational $SL(2)$-equivariant morphism $f^\pm$ is the  
contraction of the unique 
$1$-dimensional $SL(2)$-orbits $C^{\pm} \subset  E_{h,m}^\pm$ 
$( C^{\pm} \cong \P^1)$ 
into the unique $SL(2)$-fixed 
singular point $O \in E_{h,m}$.  
We remark that the $\Q$-factorial $SL(2)$-variety  $E_{h,m}^+ \cong 
SL(2) \times_B S^+ $ 
was first constructed and investigated by Panyushev in \cite{Pa88,Pa91}.
\medskip

In Section 4 we consider $SL(2)$-equivariant flips from the 
point of view of the theory of spherical varieties developed  by 
Luna and Vust \cite{LV83} (see also \cite{K91}).  
It is easy to see that any affine $SL(2)$-variety $E_{h,m}$ 
admits an additional $\C^*$-action 
which commutes with the $SL(2)$-action. If $H \subset 
SL(2) \times \C^*$ is the stabilizer subgroup of a generic 
point $x \in E_{h,m}$, then   
$(SL(2) \times \C^*)/H$ is a 
spherical homogeneous space and  $E_{h,m}$ is a 
spherical embedding. We describe simple spherical 
varietes $E_{h,m}$, $E_{h,m}^-$ and 
$E_{h,m}^+$  by colored cones. Thus, any $SL(2)$-equivariant flip 
provides an  illustration of  general results  
on Mori theory for  spherical varieties due to Brion \cite{B93},  
Brion-Knop \cite{BK94}.   
According to Alexeev and Brion \cite{AB04}, any spherical variety $V$ 
admits a flat degeneration to a toric variety $V'$. 
We apply this fact to spherical varieties $E_{h,m}$, $E_{h,m}^-$, and 
 $E_{h,m}^+$ and investigate 
the corresponding degenerations of $SL(2)$-equivariant 
flips to toric flips. We remark that the idea of toric degenerations 
appeared already in earlier papers of Popov \cite{P87} 
and Vinberg \cite{V95}. Toric degenerations of   
affine spherical varieties (including $SL(2)$-varieties) 
was considered by Arzhantsev in 
\cite{Ar99}.

\bigskip
{\bf Acknowledgments.} The authors are grateful to Ivan Arzhantsev, 
Michel Brion  and  
J\"urgen Hausen for useful discussions  and their help.

\section{Affine $SL(2)$-varieties as categorical quotients}

The complete classification of 
 normal  affine quasihomogeneous $SL(2)$-varieties has been obtained by 
Popov \cite{P73}. A shorter modern presentation 
of this classification is contained in 
the book of Kraft \cite[III.4]{Kr84}:

\begin{theorem}\cite{P73} 
Every $3$-dimensional normal affine  quasihomogeneous 
$SL(2)$-variety containing more than $1$ orbit is uniquely 
determined by a pair of numbers $(h,m) \in \{\Q \cap (0,1] \}\times \N$.  
We denote the corresponding 
variety by $E_{h,m}$.  The number $h$ is called height of  $E_{h,m}$.
The number  $m$ is called degree of  $E_{h,m}$ and it 
equals  the order of the stabilizer
of a point in the open dense $SL(2)$-orbit  ${\mathcal U} 
\subset E_{h,m}$ $($this 
stabilizer is always a cyclic group$)$. 
\label{class}
\end{theorem} 

Let $\mu_n = \langle \zeta_n \rangle$ 
be the cyclic group of $n$-th roots of unity. We denote a cyclic 
group of order $n$ also by $C_n$. We use the following notations 
for some closed subgroups in $SL(2)$: 
\[ T:= \left\{ \begin{pmatrix} t & 0 \\ 0 & t^{-1} \end{pmatrix}\;  : \; 
t \in 
\C^* \right\}, \;\; B:= \left\{ \begin{pmatrix} t & u \\ 
0 & t^{-1} \end{pmatrix}\;  : \; 
t \in 
\C^*, u \in \C \right\},\]
\[ U_n := \left\{ \begin{pmatrix} \xi & u \\ 0 & \xi^{-1} \end{pmatrix} \; : 
\; u \in 
\C, \; \xi^n =1 \right\}, \;\; 
 U:= \left\{ \begin{pmatrix} 1 & 0 \\ u & 1 \end{pmatrix} \; : 
\; u \in 
\C\right\}. \]

\begin{remark}
{\rm 
If $h =1$, then $E_{1,m}$ is smooth and it contains two 
$SL(2)$-orbits: 
$${\mathcal U} \cong SL(2)/C_m \; \mbox{\rm and} \; {\mathcal D} \cong 
SL(2)/T.$$ 
The geometric description of $E_{1,m}$ 
is easy and well-known \cite[III, 4.5]{Kr84}: 
\[ E_{1,m} \cong SL(2) \times_T  \C, \]
where the torus $T$ acts on $\C$ by character 
$\chi_m\, :\, t \to t^m$. So  $E_{1,m}$ can be considered as a 
line bundle over $SL(2)/T$. 
}
\label{h=1} 
\end{remark}

\begin{remark} 
{\rm If $0 <h<1$, then  $E_{h,m}$ contains a unique $SL(2)$-invariant 
singular point $O$.  We write $h = p/q$ 
where $g.c.d.(p,q)=1$ and define $k:=g.c.d.(m, q-p)$, $ a:= {m}/{k}$.
Then  $E_{h,m}$ contains three
$SL(2)$-orbits: 
\[  {\mathcal U} \cong 
SL(2)/C_m, \;  {\mathcal D} \cong 
SL(2)/U_{a(p+q)}, \;  \mbox{\rm and} \; \{O\}.\] 
The explicit construction of  $E_{h,m}$ given in \cite{P73} and 
\cite{Kr84} 
involves finding 
a system of generators of  the following semigroup
\[ M_{h,m}^+ := \{ (i,j) \in \Z_{\geq 0}^2 \; : \; j \leq hi, \; m | (i-j) 
\}. \]
Let $V_n$ be the standard $(n+1)$-dimensional irreducible representation 
of $SL(2)$ in the space of binary forms of degree $n$. Denote by 
$(i_1,j_1), \ldots, (i_r,j_r)$ a system of generators of 
the semigroup 
$M_{h,m}^+$. Then $E_{h,m}$ is isomorphic to  the closure $\overline{SL(2)v}$ 
of the $SL(2)$-orbit of  the vector 
\[ v := (X^{i_1}Y^{j_1}, \ldots,  X^{i_r}Y^{j_r}) \in V_{i_1+j_1} \oplus 
\cdots \oplus V_{i_r+j_r}. \]
For example, if  $m = (q-p)a, \;  a \in \N$,  then the semigroup 
 $M_{h,m}^+$ is minimally generated by $ap+1$ elements 
$(m,0), (m+1,1), (m+2,2), \ldots, (aq,ap)$ and 
$$
 v := (X^{m}, X^{m+1}Y, \ldots,  X^{aq}Y^{ap}) \in V_{m} \oplus V_{m+2}  
\oplus
\cdots \oplus V_{aq + ap} \cong V_{aq} \otimes  V_{ap}.$$
} 
\label{h<1}
\end{remark}

\begin{figure}[htbp]\label{bild1}
\begin{center}
\includegraphics[width=0.7\textwidth]{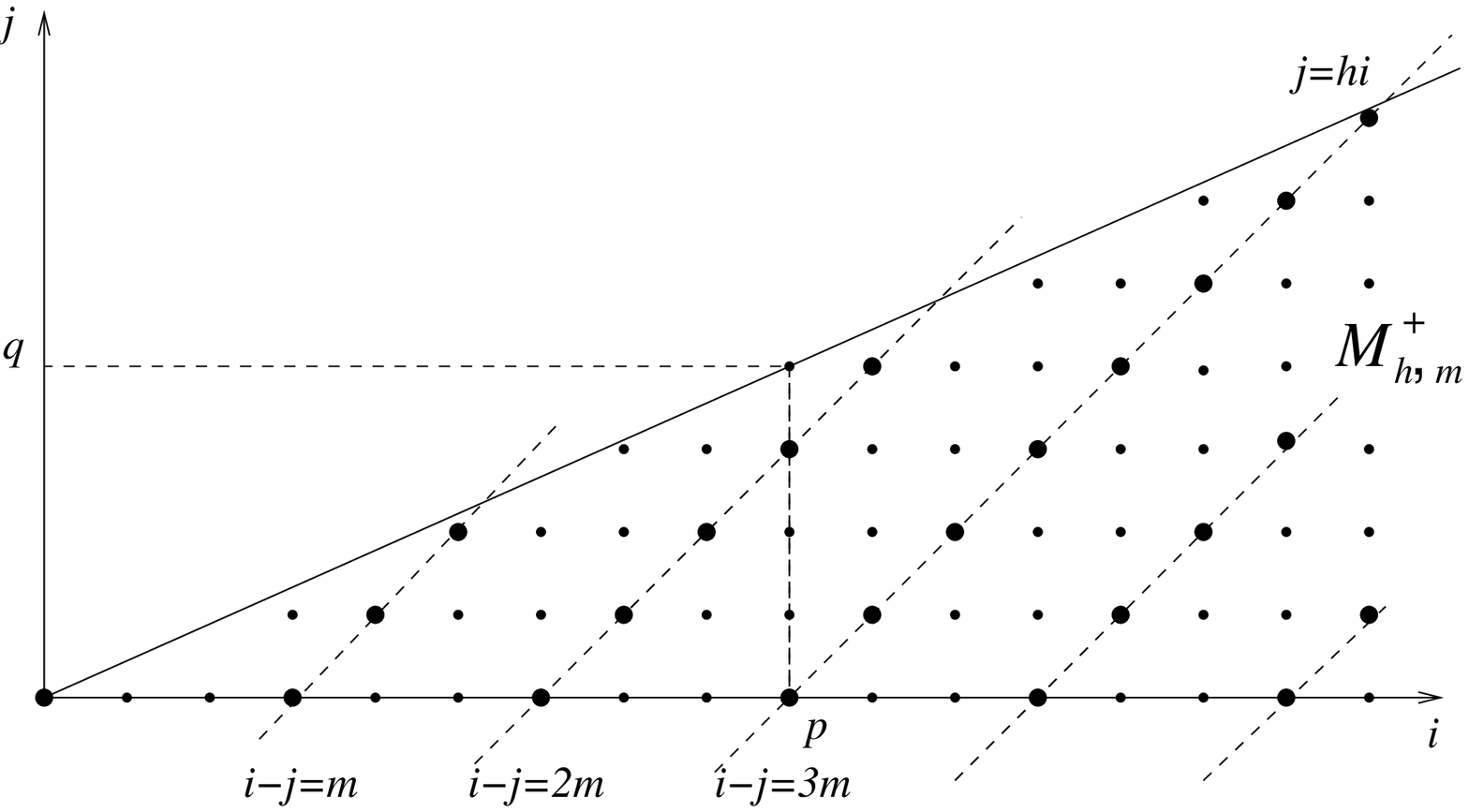}
\end{center}
\caption{}
\end{figure}

\begin{remark} \label{unique}
{\rm It is easy to see that the numbers $h$ and $m$ are uniquely 
determined by the embedding of the monoid  $M_{h,m}^+$ into  $\Z_{\geq 0}^2$ (see
Figure 1). 
} 
\end{remark}

There exists another relation between  
the submonoid  $M_{h,m}^+ \subset  \Z_{\geq 0}^2$ and $E_{h,m}$:

\begin{theorem} \cite[III,4.3]{Kr84}
Let $E$ be a normal affine $3$-dimensional quasihomogeneous  
$SL(2)$-variety with the affine coordinate ring $\C[E]$. Denote by 
$\C[E]^U$ the $U$-invariant subring. We can consider $\C[E]^U$ as 
a subring of $\C[SL(2)]^U \cong \C [X,Y]$, where $\C [X,Y]$ is 
the algebra of regular functions on $SL(2)/U \cong \C^2 \setminus 
\{(0,0)\}$. Then the monomials $\{ X^iY^j\; | \; (i,j) \in M_{h,m}^+ \}$
form a $\C$-basis of $\C[E]^U$.
\label{U-inv}  
\end{theorem} 
 
\medskip

Our next purpose is to give a new description of an  
affine quasihomogeneous  
$SL(2)$-variety $E_{h,m}$ as a categorial quotient of a $4$-dimensional 
affine hypersurface. Throughout this paper by a categorical quotient 
$X/\!\!/G$ 
of an irreducible affine algebraic variety $X$ over $\C$ by a reductive group 
$G$ we mean ${\rm Spec}\, \C[X]^G$ \cite{FM82,VP89}.   

Let $SL(2) \times \C^5 \to \C^5$ be the  
$SL(2)$-action on $\C^5$ considered as 
$V_0 \oplus V_1 \oplus V_1$. We use the coordinates $X_0, X_1,X_2,X_3, X_4$
on $\C^5$ and identify $X_1, X_2, X_3, X_4$ with the coefficients 
of the $2\times 2$-matrix
\[ \begin{pmatrix} X_1 & X_3 \\ X_2 & X_4 \end{pmatrix} \]
on which $SL(2)$ acts by left multiplication. Denote by $D(5, \C)$ 
the group of diagonal matrices of order $5$ acting on $\C^5$.

One of our main results of this section is the following: 

\begin{theorem} 
Let $E_{h,m}$ be a normal  affine $SL(2)$-variety of height $h = p/q \leq 1$ 
$(g.c.d.(p,q)=1)$ and of 
degree $m$. 
Then $E_{h,m}$ is isomorphic 
to the categorical quotient of the affine hypersurface $H_{q-p} \subset \C^5$ 
defined by the equation 
\[ X^{q-p}_0 = X_1 X_4 - X_2X_3 \]
modulo the action of the diagonalizable group 
$G_0 \times G_m \subset D(5,\C)$, 
where $G_0 \cong \C^*$ 
consists of diagonal matrices 
 $\{ diag(t,t^{-p},t^{-p},t^q,t^q)\; : \; t \in \C^* \}$ 
and $G_m \cong \mu_m = 
\langle \zeta_m \rangle$ is generated by 
$diag(1,\zeta_m^{-1}, \zeta_m^{-1}, \zeta_m, \zeta_m )$. 
\label{cat1}
\end{theorem}

\noindent
{\em Proof.} 
{\sc Case 1}: $h =1$. Then $p=q =1$,  
\[ G_0:= \{ diag (t,t^{-1},t^{-1}, t,t) \; : \; t \in \C^* \}, 
\;\;  G_m:= \{ diag (1,\zeta^{-1},\zeta^{-1}, \zeta, \zeta) 
\; : \; \zeta  
\in \mu_m \}, \]
and the hypersurface $H_0$ is defined by the equation 
$1 = X_1 X_4 - X_2X_3$. 
The  algebraic group $G_0 \times G_m$ 
can be written as a direct product in another way:
\[ G_0 \times G_m  = G_0 \times G_m', \]
where $G_m'  := \{ diag(\zeta, 1,1,1,1)\; : \;\zeta \in \mu_m  \}$.  
 We remark that the hypersurface $H_0$ is isomorphic 
to the product $SL(2) \times \C$. Moreover,  the $G_0$-action on the 
first factor 
$SL(2)$ is the same as the action of the maximal torus $T$ 
by  right 
multiplication. On the other hand, $H_0 /\!\!/ {G_m'}$ is again 
isomorphic to $SL(2) \times \C$, because 
$G_m'$ acts trivially on $SL(2)$ and $\C/\!\!/ {G_m'} \cong \C$ 
(one replaces the coordinate $X_0$ on $\C$ by a new 
$G_m'$-invariant coordinate $Y_0 = X_0^m$). 
So the $G_0$-action 
on the second factor $\C$ in $SL(2) 
\times \C \cong H_0/\!\!/{G_m'}$ is defined  by 
the character $\chi_m\; : \; t \to t^m$ .  
Thus,  we come to the already known 
description of $E_{1,m}$ as a $T$-quotient: $E_{1,m} 
\cong SL(2)\times_T \C$ (see \ref{h=1}).

{\sc Case 2}: $m=1$, $h =p/q <1$. 
The $SL(2)$-action on $\C ^5$  commutes with the $G_0$-action 
and  the hypersuface $H_{q-p}$ defined by the equation 
\[ X^{q-p}_0 = X_1 X_4 - X_2X_3. \]
is invariant under this $G_0 \times SL(2)$-action. 
Moreover, the point $x :=(1,1,0,0,1) \in H_{q-p}$ has a trivial stabilizer
in  $G_0 \times SL(2)$. Therefore $H_{q-p}$ is the  closure 
of  the $G_0 \times SL(2)$-orbit of $x$ in $\C^5$ and $X_{p,q}:= 
{\rm Spec}\, \C[H_{q-p}]^{G_0}$ is an affine $SL(2)$-embedding. 
One can identify 
the open dense $SL(2)$-orbit ${\mathcal U}$ in $X_{p,q}$ with the 
$G_0$-quotient 
of the open subset in $H_{q-p}$ defined by the 
condition  $X_0 \neq 0$. Moreover, the affine coordinate ring 
$\C[{\mathcal U}]$ is generated by the $G_0$-invariant monomials
\begin{equation} X:= X_0^pX_1,\;\; Y:= X_0^{-q}X_3,\;\; Z:= X_0^pX_2,\;\;
 W:= X_0^{-q}X_4 \end{equation}
satisfying the equation 
\[ \det   \begin{pmatrix} X & Y \\ Z & W \end{pmatrix} = X_0^{p-q}X_1X_4 - 
X_0^{p-q}X_2X_3 = 1. \]
By a theorem 
of Luna-Vust \cite[III,3.3]{Kr84}, the normality of  $X_{p,q}:= 
{\rm Spec}\, \C[H_{q-p}]^{G_0}$ follows
from the normality of  $\C[X_{p,q}]^U = 
\C[ H_{q-p}]^{G_0 \times U}$. It is easy to see that 
$$\C[H_{q-p}]^U \cong  \C[X_0,X_1,X_3].$$ Since $U$-action and 
$G_0$-action commute, it remains to compute 
the $G_0$-invariant subring  $\C[X_0,X_1,X_3]^{G_0}$
under the 
$\C^*$-action of $G_0$ on $\C^3$ defined by $diag(t,t^{-p}, t^q)$. 
Straightforward calculations show that 
the ring  $\C[X_0,X_1,X_3]^{G_0}$ has a $\C$-basis 
consisting of all monomials $X^iY^j = 
X_0^{pi-qj}X_1^iX_3^j \in  \C[{\mathcal U}]^U = \C[X,Y]$ 
such that $pi -qj \geq 0$, $i \geq 0$, $j \geq 0$, i.e. , 
$(i,j) \in 
M_{h,1}^+$. By \ref{unique} and \ref{U-inv}, we obtain 
 simultaniously that $X_{p,q}$ is normal and 
that $X_{p,q} \cong E_{h,1}$. 

{\sc Case 3}: $m>1$, $h =p/q <1$. Let $X_{p,q}^m$ 
be the categorical quotient of $H_{q-p}$ by $G_0 \times G_m$ where 
  $G_m \cong \mu_m = \langle 
\zeta_m \rangle$ acts   
by 
$diag(1, \zeta_m^{-1},\zeta_m^{-1},\zeta_m,\zeta_m)$. 
By the same arguments as above, 
one obtains that $\C[X_{p,q}^m]^U \cong  
\C[X_0,X_1,X_3]^{G_0 \times G_m}$ where $G_m \cong \mu_m$ acts on $\C^3$ 
by $diag(1,\zeta_m^{-1},\zeta_m)$ and $G_0$ on $\C^3$ 
by $diag(t,t^{-p},t^q)$. Therefore the ring $\C[X_{p,q}^m]^U \subset 
\C[{\mathcal U}]^U = \C[X,Y]$ has 
a $\C$-basis 
consisting of all monomials 
$X^iY^j = X_0^{pi-qj}X_1^iX_3^j$ such that $(i,j) \in 
M_{h,m}^+$ (the condition  $m | (j-i)$ follows from the  
$G_m$-invariance of monomials  $X_0^{pi-qj}X_1^iX_3^j$). By  \ref{unique} 
and \ref{U-inv}, this 
shows that  $X_{p,q}^m \cong E_{h,m}$. 
\hfill $\Box$ 
\medskip

It will be important to have the following another similar description of an 
arbitrary affine normal quasihomogeneous $SL(2)$-variety  $E_{h,m}$ 
as a categorical quotient of an affine  hypersurface:

\begin{theorem} 
Let $E_{h,m}$ be a normal  affine $SL(2)$-variety of height $h = p/q \leq 1$ 
$(g.c.d.(p,q)=1)$ and of 
degree $m$. We define $ b:= (q-p)/k$. 
Then $E_{h,m}$ is isomorphic 
to the categorical quotient of the affine hypersurface $H_b \subset \C^5$ 
defined by the equation 
\[ Y_0^b = X_1 X_4 - X_2X_3 \]
modulo the action of the diagonalizable group  
$G:= G_0' \times G_a \subset D(5,\C)$, 
where $G_0' \cong \C^*$ 
consists of diagonal matrices 
 $\{ diag(t^k,t^{-p},t^{-p},t^q,t^q)\; : \; t \in \C^* \}$ 
and $G_a \cong \mu_a = 
\langle \zeta_a \rangle$ is generated by 
$diag(1,\zeta_a^{-1}, \zeta_a^{-1}, \zeta_a, \zeta_a )$. 
\label{cat2}
\end{theorem}

\noindent
{\em Proof.} 
By \ref{cat1}, we have $E_{h,m} = H_{q-p} /\!\!/(G_0 \times G_m)$. 
We note that the conditions  
$k:=g.c.d.(q-p,m)$ and $g.c.d.(q,p)=1$ imply that $g.c.d.(k,p) = 
g.c.d.(k,q)=1$. Since $\zeta_m^a$ is a generator of $\mu_k$ and since
the maps $z \to z^p$ and $z \to z^q$ are bijective on  $\mu_k$ we can find 
another  generator $\xi \in \mu_k$ such that 
$\xi^p \zeta_m^a =  \xi^q \zeta_m^a =1$. 
Therefore,  
$G_0 \times G_m$ contains the following element
\[ g =  diag(\xi, 1,1,1,1) = (\xi, \xi^{-p}, \xi^{-p},\xi^{q},\xi^{q}) 
\cdot (1,\zeta_m^{-a},\zeta_m^{-a}, \zeta_m^{a},\zeta_m^{a}).\]
Consider the homomorphism 
\[ \psi_k\; : \; D(5, \C) \to D(5, \C), \; \; 
(\lambda_0, \lambda_1, \lambda_2, \lambda_3, \lambda_4) \mapsto 
(\lambda_0^k, \lambda_1, \lambda_2, \lambda_3, \lambda_4). \]
Then $\psi_k(G_0) = G_0'$ and 
$$ G'_k:= {\rm Ker}\, \psi_k \cap (G_0 \times G_m) = 
\langle g \rangle= 
\{ diag(\zeta, 1,1,1,1)\; : \; \zeta \in \mu_k \}.$$
So we obtain a short exact sequence
\[  1 \to G_k' \to G_0 \times G_m \to  G_0' \times G_a \to 1, \]
where 
$$G_a = \{ diag(1,\zeta^{-1},\zeta^{-1},\zeta, \zeta)\; 
:\; \zeta \in \mu_a \}.$$
Therefore the categorical  
$G$-quotient of $H_{q-p}$ can be divided 
in two steps. First we divide  $H_{q-p}$ by  the subgroup $G_k' \subset 
 G_0 \times G_m$  
and after that divide   by 
the group $G_0' \times G_a$. Using a  new $G_k'$-invariant 
coordinate $Y_0 = X_0^k$, we see that  $H_{q-p} /\!\!/{G_k'}$ is 
isomorphic to the hypersurface $H_{b}$  
defined by the equation 
\[ Y_0^b = X_1X_4 - X_2X_3. \]
Since $G_0$ acts on $Y_0$ by character $t \to t^k$, 
$E_{h,m} \cong X_{p,q}^m =  H_{q-p} /\!\!/(G_0 \times G_m)$ 
is isomorphic to the categorical quotient of $H_b$ 
modulo   the above $G_0' \times G_a$-action.  
\hfill $\Box$

\section{The Cox ring of an affine $SL(2)$-variety}

Let us review a  definition of the total coordinate ring (or Cox ring) 
of a normal algebraic variety $X$ with finitely generated 
divisor class group ${\rm Cl}(X)$ (see e.g. \cite{Ar08,H08}). 

\begin{definition} 
{\rm Let $X$ be a normal quasiprojective irreducible 
algebraic variety over $\C$ with the field of rational functions 
$\C(X)$. 
We assume 
that ${\rm Cl}(X)$ is a finitely generated abelian group and that every 
invertible regular function on $X$ is constant.  
Choose divisors  $D_1, \ldots, D_r$  in $X$ whose 
classes generate  ${\rm Cl}(X)$ and  prinicipal divisors $D_1',\ldots , 
D_s'$ which     
form a $\Z$-basis of the kernel 
of the surjective homomorphism $\varphi\, : \, \Z^{r} 
\to {\rm Cl}(X)$. Furthermore, 
we choose rational functions $f_1, \ldots, f_s \in 
\C(X)$ such that $D_i' = (f_i)$   $(i =1, \ldots, s)$. 
For any $k = (k_1, \ldots, k_r ) \in \Z^r$, we consider 
a divisor $D(k) := \sum_{j =1}^r k_j D_j$ 
and put 
$${\mathcal L}(D(k)) := \{ f \in \C(X)\; : \; D(k) + (f) \geq 0.\}$$ 
Then for every $i \in 
\{ 1, \ldots, s \}$, one has  an isomorphism 
\[ \alpha_i \; : \;     
{\mathcal L}(D(k))  \stackrel{\cong}{\to} {\mathcal L}(D(k)+D_i')
, \;  \alpha_{i}( f ) = \frac{f}{f_i} \;\; \forall f \in 
{\mathcal L}(D(k)).  \]    
In the $\Z^{r}$-graded ring 
\[  {\mathcal R}:= \bigoplus_{k \in {\Z}^r }   {\mathcal L}(D(k)), \]
we consider the  ideal ${\mathcal I}$ generated by  
all elements $f -  \alpha_{i}(f)$ $\forall f 
\in {\mathcal L}(D(k)) $, $\forall k \in {\Z}^r$, and 
 $\forall i \in \{ 1, \ldots, s\}$.  
The ring $${\rm Cox}(X) := {\mathcal R}/{\mathcal I}$$ is called 
{\bf Cox ring} of $X$ associated with divisors  
$D_1, \ldots, D_r$ and rational functions 
$f_1, \ldots, f_s$. By \cite[Prop. 3.2]{Ar08},   ${\rm Cox}(X)$ is 
uniquely defined up to isomorphism and does not 
depend  on the choice of  
generators $D_1, \ldots, D_r$  
of  ${\rm Cl}(X)$ and rational functions  $f_1, \ldots, f_s$. 
Moreover, one has a natural 
${\rm Cl}(X)$-grading   
\[ {\rm Cox}(X) =  {\mathcal R}/{\mathcal I} 
\cong   \bigoplus_{c \in {\rm Cl}(X)} 
\Gamma(X, {\mathcal O}_X(c)). \]
} 
\label{cox-def}
\end{definition} 
\medskip

Let $A$ be a finitely generated abelian group. 
We shall need the following  criterion  for a finitely generated factorial 
$A$-graded $\C$-algebra $R$ with $R^\times = \C^*$ to be a Cox ring of a normal 
quasiprojective algebraic variety $X$ with $A \cong {\rm Cl}(X)$.

\begin{theorem} 
Let $Y$ be a normal irreducible affine algebraic variety over $\C$ with 
a factorial coordinate ring $R=\C[Y]$. We assume that  
$\Gamma (Y, {\mathcal O}_Y^*) = \C^*$ 
and that $Y$ admits a regular action  $G \times Y \to Y$ 
of a diagonalizable group $G$, or, equivalently, $R$ admits 
an $A$-grading by the group $A = {\rm Hom}_{\rm alg}(G, \C^*)$ 
of algebraic characters of $G$.  
Then $R$ is a Cox ring 
of some normal quasiprojective algebraic variety $X$  such that 
${\rm Cl}(X) \cong  A$ and $\Gamma (X, {\mathcal O}_X^*) = 
\C^*$  if and only the following 
conditions are satisfied:  

(i) there exists  an open dense nonsingular $G$-invariant subset 
$U \subset Y$ such that   $\codim_Y 
(Y \setminus U) \geq 2$ and  $G$ acts freely on $U$; 

(ii) there exists a character $\chi \in {\rm Hom}_{\rm alg}(G, \C^*)$ such that 
$U \subset Y^{ss}(L)$, where $L$ is the $G$-linearization of 
the trivial line bundle over $Y$ corresponding to $\chi$. 
\label{criterion}
\end{theorem} 

\noindent
{\em Proof.} Assume that $Y$ admits a regular $G$-action such that 
the  conditions (i), (ii) are satisfied. We define $X$ to be 
$Y^{ss}(L)/\!\!/G$. Then $X$ is a normal 
irreducible quasiprojective variety and $\Gamma (X, {\mathcal O}_X^*) 
= \Gamma (Y, {\mathcal O}_Y^*)^G = 
\C^*$. Moreover, $\overline{U}:= U/G$ is a smooth open subset of $X$.     
Let us show that ${\rm Cl}(X) \cong  A$, where  
$A = {\rm Hom}_{\rm alg}(G, \C^*)$.  

Since $R$ is factorial and $U$ is a smooth open subset of $Y$, 
we have ${\rm Pic}(U) = {\rm Cl}(U)=0$. By a general result 
in \cite[5.1]{KKV89}, the Picard group of  $\overline{U}$ 
is isomorphic to the group of $G$-linearizations of the trivial line  
bundle over $U$. On the other hand, since ${\rm codim}_Y\, (Y \setminus U) 
\geq 2$ and $Y$ is normal, all invertible 
regular functions on $U$ extend to invertible regular 
functions on $Y$, i.e., they are 
constant. By  \cite{KKV89}, the latter implies 
that the group of $G$-linearizations of the trivial line bundle over 
$U$ is isomorphic to the group of characters of $G$, i.e., 
\[ {\rm Pic}(\overline{U}) \cong 
{\rm Hom}_{\rm alg}(G, \C^*)  = A. \]
Since ${\rm Pic}(\overline{U}) = {\rm Cl}(\overline{U} )$, 
it remains to show that 
${\rm codim}_X (X \setminus \overline{U}) \geq 2$. Assume
 that there exists an 
irreducible nonempty divisor $Z \subset X$ such that   
$ \overline{U}\cap Z = \emptyset$. 
Since $X$ is normal, the local ring ${\mathcal O}_{X,Z}$ is a discrete 
valuation ring, i.e., there exists an affine open subset $U' \subset X$
such that $Z' := U' \cap Z \neq \emptyset$, $\overline{U} \cap Z' = 
\emptyset$,  
and $Z'$ is a principle divisor in 
$U'$ defined by a regular fuction $g \in \C[U']$. 
Consider the morphism 
\[ \pi \; : \; Y^{ss}(L) \to X. \] 
Without loss of generality, we can assume 
$\widetilde{U'}:= \pi^{-1}(U')$ is an affine open subset in  $Y^{ss}(L)$ 
and $\C[U'] = \C[ \widetilde{U'}]^G$. Then the element 
$\tilde{g}:= \pi^*(g) \in  \C[ \widetilde{U'}]$ defines a principle 
divisor $\widetilde{Z'} := (\tilde{g})   
\subset \widetilde{U'}$ such that 
$\widetilde{Z'}  \cap U = \emptyset$ and $\widetilde{Z'} \neq 
\emptyset$. The latter contradicts to 
 $\codim_{\widetilde{U'}} 
(\widetilde{U'} \setminus (\widetilde{U'}\cap U)) \geq 
\codim_Y 
(Y \setminus U) \geq 2$, i.e., we must have $Z' =\emptyset$. 

In order to identify $R = \bigoplus_{a \in A} R_a$  
with  the Cox ring of $X$ we consider a finite subset  
$\{ a_1, \ldots, a_r \} \subset A$ such that the homogeneous 
components $R_{a_1}, \ldots, R_{a_r}$ generate the algebra $R$ and 
$R_{a_i} \neq 0$ 
for all $i \in \{1, \ldots, r \}$. 
Since the class of any effective divisor in $X$ is a nonnegative 
integral linear combination of  $a_1, \ldots, a_r$, we obtain that 
$a_1, \ldots, a_r$ are generators of $A$. 
We choose $r$ nonzero 
elements $g_j \in R_{a_j}$, $j \in \{ 1, \ldots, r \}$ 
which define 
$r$ effective principal divisors $\widetilde{D_j} = (g_j)$ 
in $Y$ 
 ($j \in 
\{ 1, \ldots, r\})$.
Then  we obtain $r$ effective divisors in $X$: 
\[ D_j :=  (\widetilde{D_j} \cap Y^{ss}(L))/\!\!/G ,\;\;  
j  \in \{ 1, \ldots, r \}. \]
Consider the epimorphism 
$\varphi \; : \; \Z^{r} \to A$. 
For any $k = (k_1, \ldots, k_r) \in \Z^{r}$ we define a rational function
$$g(k):= g_{1}^{k_1} \cdots g_{r}^{k_r}  \in \C(Y)$$
and a divisor 
$$  D(k) := k_1 D_{1} + \cdots + k_r D_{r} 
\in {\rm Div}(X).$$
If $a_1', \ldots, a_s'$ is a $\Z$-basis of ${\rm Ker}\, \varphi$, then 
$s$ rational functions $f_i:= g(a_i')$ $(i =1, \ldots, s)$ are $G$-invariant, 
i.e., elements of  $\C(X)$. 
So we obtain $s$ principle divisors $D_i':= D(a_i') = (f_i)$ in $X$.
On the other hand, for any  $k \in \Z^{r}$, one has
\[ {\mathcal L}(D(k)) =\left\{ \frac{h}{g({k})} \in \C(X) \; : \; 
h \in R_{\varphi(k)} \right\}. \]
Consider the $\Z^r$-graded ring  
\[ {\mathcal R}:= \bigoplus_{k \in  \Z^{r}}  {\mathcal L}(D(k)) \]
together with the surjective homogeneous homomorphism 
\[ \beta\; : \;  {\mathcal R} \to R = \bigoplus_{a \in A} R_a \]
whose restriction to $k$-th homogeneous component is an isomorphism 
\[ \beta_k \; : \;  {\mathcal L}(D(k)) \stackrel{\cong} 
\to  R_{\varphi(k)} \]
defined by multiplication with $g(k)$. Then the elements 
\[ \left( \frac{h}{g(k)} - \frac{h}{g(k +a_i')}\right) = 
\left( \frac{h}{g(k)} - \frac{h}{g(k)f_i}\right) 
\in 
{\mathcal R}_{k} \oplus {\mathcal R}_{k +a_i'}
, \]
\[ \forall k \in \Z^r, \; 
\forall h \in R_{\varphi(k)} = R_{\varphi(k+a_i')}, \;  \forall i \in 
\{ 1, \ldots, s \}\] 
are contained in ${\rm Ker}\, \beta$. Therefore, $\beta$ induces 
a surjective homogeneous homomorphism of the Cox ring 
${\mathcal R}/{\mathcal I}$ to $R$. By comparing the homogeneous 
components of ${\mathcal R}/{\mathcal I}$ and  $R$, we obtain an isomorphism
${\mathcal R}/{\mathcal I} \cong R$. 
\medskip

Now assume that a factorial  $A$-graded 
$\C$-algebra $R$ is the Cox ring of some 
normal irreducible quasiprojective variety $X$ with ${\rm Cl}(X) \cong A$. 
Using the same idea as in \ref{cox-def}, we can define  
a sheaf-theoretical version of the Cox ring of $X$ (see 
\cite[Section 2]{H08}):
 $$\widetilde{R} = \bigoplus_{a \in A} {\mathcal O}_X(a) $$ 
which is a $A$-graded ${\mathcal O}_X$-algebra such that 
$\Gamma(X, \widetilde{R}) = R$. 
Define $Y':= {\rm Spec}_X (  \widetilde{R})$ as a relative spectrum over 
$X$. By \cite[Prop.2.2]{H08},  
$Y' \subset Y:= {\rm Spec}_{\C}(R)$ is an open embedding and  the morphism 
$\pi\; : \; Y' \to X$ is a categorical quotient by the action 
of $G:= {\rm Spec}\, \C[A]$. Moreover, $G$ acts freely on the open 
subset $U:= \pi^{-1}(\overline{U})$, where  
$\overline{U}:= X \setminus {\rm Sing}(X) \subset X$
 the set of all smooth points of $X$ and ${\rm codim}_Y  
(Y \setminus U) \geq 2$. 
We consider a locally closed
embedding $\jmath\, : \, X \to \P^n$  
and define  ${\mathcal L}:= \jmath^* {\mathcal O}(1)$.  
Since ${\rm Cl}(Y') =0$, the pullback  $L:= \pi^*{\mathcal L}$ is 
a trivial line bundle over $Y'$ having a $G$-linearization. Since 
all invertible global regular functions on $Y'$ are constants, 
this $G$-linearization is determined by a character $\chi \in 
{\rm Hom}_{\rm alg} (G, \C^*) \cong  A$. Since $\pi \, : \, 
Y' \to X$ is a categorical quotient, we have $U \subseteq  Y'\subseteq 
Y^{ss}(L)$. Theorem is proved. 
\hfill $\Box$ 
\bigskip

\begin{remark} 
{\rm Methods in \cite{H08} allow to formulate and prove a more 
general version of \ref{criterion} for algebraic varieties $X$ which 
are not necessary  quasiprojective. Moreover, in Theorem  
\ref{criterion} it is enough to assume 
only  $A$-graded factoriality  
of $R$, i.e., that  every $A$-homogeneous divisorial ideal is principal. 
} 
\end{remark}

Now we  begin with the following observation:

\begin{proposition}
The affine coordinate ring $\C[H_b]$ of the hypersurface $H_b \subset \C^5$ 
is factorial. Invertible elements in  $\C[H_b]$ are exactly nonzero constants.
\label{factor}
\end{proposition} 

\noindent
{\em Proof.} 
Consider the open subset  $U_2^+ \subset H_b$ 
  defined by  $X_2 \neq 0$. Since  $U_2^+$ is 
isomorphic to a Zariski  open subset 
in $\C^4$, we obtain  $ {\rm Cl}(U_2^+)  =0$. The complement  
$\widetilde{S^+} := H_b \setminus U_2^+$ is a principle divisor $(X_2)$. 
We note that $\widetilde{S^+}$ defined by the binomial equation 
$Y_0^b = X_1 X_4$ which shows that $\widetilde{S^+}$ is isomorphic 
to the product 
of $\C$ (with the coordinate $X_3$) and a 
$2$-dimensional affine toric variety with 
a $A_{b-1}$-singularity defined by the equation $Y_0^b = X_1 X_4$.  
Therefore, $\widetilde{S^+}$ is irreducible 
and the short exact localization sequence

\[
\begin{matrix} \Z & \to & {\rm Cl}(H_b) & \to  & {\rm Cl}(U_2^+) & \to & 0 \\
1 &  \mapsto &  [\widetilde{S^+}]  & & & &
\end{matrix}
\] 
shows that $[\widetilde{S^+}] = 0 \in {\rm Cl}(H_b)$, i.e., 
  the image of $\Z$ in  $ {\rm Cl} (H_b)$ is zero. 
Thus, we obtain $ {\rm Cl}(H_b) =0$. 

In order to prove the second statement we consider the following two cases. 

Case 1: $b =0$. Then $H_b \cong SL(2) \times \C$. Since all invertible 
elements in the coordinate ring of $SL(2)$ are constants, we obtain 
the same property for  the coordinate ring of $SL(2)\times \C$. 

Case 2: $b > 0$. Then we can define a $\Z_{\geq 0}$-grading 
of $\C[H_b]$ by setting 
$\deg\, X_1 = \deg\, X_2 = \deg\, X_3 = \deg \,X_4 = b$ and $\deg\, Y_0 = 2$. 
Since the $0$-degree component of  $\C[H_b]$ is $\C$, we obtain 
that all invertible elements in $\C[H_b]$ are nonzero constants. 
\hfill $\Box$
\medskip

\begin{proposition} 
Consider   the following two Zariski open subsets $U^+$, $U^-$ 
in $H_b$:
\[ U^+ := H_b \setminus \{X_1 =X_2 =0\}, \;\; \;  U^- := H_b \setminus 
\{X_3 =X_4 =0\}. \]
Denote   $U_{h,m} :=  E_{h,m} \setminus 
{\rm Sing}( E_{h,m}) \subseteq E_{h,m}$, where ${\rm Sing}( E_{h,m}) = 
\emptyset$ if $h =1$ and  ${\rm Sing}( E_{h,m}) = 
\{ O \}$ if $h <1$. 
Then  the diagonalizable 
group  $G:= G_0' \times G_a$ acts freely on $U^+ \cap U^-$ 
and  $U_{h,m} \cong (U^+ \cap U^-) /G$.
\label{free}
\end{proposition}

\noindent
{\em Proof.} Let $x = (y_0,x_1, x_2, x_3, x_4)$ be a point 
in $U^+ \cap U^-$ and $g \in G$ an element such that 
$gx =x$. We write $g$ as 
\[ g = diag(t^k, t^{-p}, t^{-p}, t^q, t^q) \cdot 
diag(1 , \zeta^{-s}, \zeta^{-s}, \zeta^{s}, \zeta^{s}), \;\; t \in
\C^*, 
\zeta \in \mu_a. \]    
Then $t^{-p}\zeta^{-s} =1$ (because at least one of $x_1$ and $x_2$ is 
nonzero), and  $t^q \zeta^{s} =1$ (because at least one of $x_3$ and $x_4$ 
is nonzero). Therefore, $t^{-p}\zeta^{-s} t^q \zeta^{s}= t^{q-p} =1$. 
Since $q-p$ and $a$ are coprime we obtain  that 
$t^{p} = \zeta^{s} = t^q =1$. Since $g.c.d.(p,q) =1$ we get 
$t =1$. Therefore $g =1$, i.e., $G$ acts freely  on $U^+ \cap U^-$.

Now we remark that the  open subsets $U^+, U^- \subset H_b$ are 
$SL(2)$-invariant and have nonempty intersection with the 
$SL(2)$-invariant 
divisor $\widetilde{D}:= 
\{ Y_0 = 0 \} \subset H_b$. Therefore, the smooth $SL(2)$-variety  
$(U^+ \cap U^-) /G$ contains more than one $SL(2)$-orbit.  So 
$(U^+ \cap U^-) /G$ coincides with 
$E_{h,m} \setminus {\rm Sing}( E_{h,m}) = U_{h,m} $ (see \ref{h=1} and 
\ref{h<1}). 
\hfill $\Box$ 
\medskip

\begin{corollary}
For any affine $SL(2)$-variety  $E_{h,m}$, one has 
\[  {\rm Cox}(E_{h,m}) \cong \C[H_b] = 
 \C[Y_0, X_1, X_2, X_3, X_4] / (  Y^b_0 - X_1 X_4 + X_2X_3). \]
\label{cox}
\end{corollary}

\noindent
{\em Proof.} 
Let $L_0$ be trivial $G$-linearized line bundle over $H_b$, i.e., 
${\mathcal O}_{H_b} \cong {\mathcal O}_{H_b}(L_0)$ as $G$-bundles. 
Then $H_b^{ss}(L_0) = H_b$ and $H_b^{ss}(L_0)/\!\!/G \cong E_{h,m}$. 
By \ref{free}, $G$ acts freely on the open subset $U:= U^+ \cap U^- \subset 
H_b$ and ${\rm codim}_{H_b}\, (H_b \setminus U) =2$. 
By \ref{criterion}, the affine coordinate ring of $H_b$ is 
isomorphic to the Cox ring of $E_{h,m}$. 
\hfill $\Box$

\begin{corollary} \cite{Ga08}
An affine $SL(2)$-variety  $E_{h,m}$ is toric if and only if $b=1$, 
i.e., $q-p$ divides $m$. 
\end{corollary} 

\noindent
{\em Proof.} If $b =0$ (i.e. $h =1$), then $E_{1,m}$ is smooth and 
 ${\rm Cl}(E_{1,m}) \cong \Z$. However, the divisor class 
group of any smooth affine toric variety is trivial. Hence,  $E_{1,m}$ 
is not toric. 

In general, if $X$ is a normal affine toric  variety such that 
all invertible elements in $\C[X]$ are constant, then ${\rm Cox}(X)$ 
is a polynomial ring \cite{Cox95}. In particular, the spectrum of 
${\rm Cox}(X)$ is nonsingular. 
On the other hand, if $b >1$, then the  
hypersurface $H_b \subset \C^5$ defined by the 
equation  $Y^b_0 - X_1 X_4 + X_2X_3 =0$ is singular. Therefore, 
$E_{h,m}$ is not toric if $b >1$. 

If $b =1$, then  $H_b \cong \C^4$, so   
$E_{h,m} \cong \C^4/\!\!/G$ is toric.
\hfill $\Box$ 
\medskip

Using \ref{cox}, we obtain a simple interpretation 
of  the following computation of ${\rm Cl}(E_{h,m})$ due 
to Panyushev:

\begin{proposition} \cite[Th.2]{Pa92} For any normal affine $SL(2)$-variety 
$E_{h,m}$, one has 
\[  {\rm Cl}(E_{h,m}) \cong  \Z \oplus C_a. \]
Let 
$D \subset  E_{h,m}$ be the closure of the unique 
$2$-dimensional $SL(2)$-orbit ${\mathcal D}$. Denote by  
$S^+ \subset E_{h,m}$ (respectively by $S^- \subset E_{h,m}$) be the closure 
  in  $E_{h,m}$ of 
the $B$-orbit in ${\mathcal U} \cong SL(2)/C_m$ defined by  
the equation $Z^m =0$  (respectively, by  $W^m =0$). Then  ${\rm Cl}(E_{h,m})$ 
is generated by two elements $[D]$ and $[S^+]$, or, respectively, by   
$[D]$ and $[S^-]$)  
satisfying the unique relation: 
\[ ap[{D}] + m [S^+] =0, \]
or, respectively,  
\[ -aq[ D] +  m [S^-] =0. \]
\label{cl-div}
\end{proposition}

\noindent
{\em Proof.} The isomorphisms  
\[ {\rm Cl}(E_{h,m}) \cong {\rm Hom}_{\rm alg}(G, \C^*) 
\cong {\rm Hom}_{\rm alg}(G_0', \C^*) \oplus 
 {\rm Hom}_{\rm alg}(G_a, \C^*) \cong \Z \oplus C_a. \]
follow immediately from \ref{cox}. 
Let $D' \subset E_{h,m}$ be an arbitrary nonzero 
effective irreducible divisor. Consider 
the surjective morphism $\pi\; : \; U^- \cap U^+ \to  (U^- \cap U^+) /G 
= U_{h,m}$. Then the support of $D'$ has a nonempty intersection 
with $U_{h,m}$, because 
${\rm codim}_{E_{h,m}}{\rm Sing}(E_{h,m}) \geq 2$. Then the closure $\widetilde{D'}$  
of $\pi^{-1}( D' \cap U_{h,m}) \subset H_b$ is a $G$-invariant principal   
irreducible divisor  (see \ref{factor}). Therefore,  $\widetilde{D'}$
is defined by zeros of a polynomial $\tilde{f}(Y_0, X_1, X_2, X_3, X_4)$ 
such that $\tilde{f}(gx) = \tilde{\chi}(g) \tilde{f}(x)$ and   
$\tilde{\chi} = \chi_{D'} 
\in {\rm Hom}_{\rm alg}(G, \C^*)$ is the character representing 
the class $[D'] \in  {\rm Cl}(E_{h,m})$. 

It is easy to see that the irreducible divisors 
$\widetilde{D}, \widetilde{S^+}, \widetilde{S^-} \subset H_b$ 
are defined respectively 
by polynomials $ Y_0, X_2, X_4 $. The corresponding 
characters $\tilde{\chi}$ of $G \cong \C^* \times \mu_a$ are : 
\[ \chi_D (t, \zeta) = t^k, \;\; \chi_{S^+}(t, \zeta) = t^{-p}\zeta^{-1}, \;\; 
\chi_{S^-}(t, \zeta)=  t^q \zeta. \]
Since $g.c.d. (ap,k) = g.c.d. (aq,k)=1$ each  pair 
$\{ \chi_D ,  \chi_{S^+} \}$ and  $\{ \chi_D ,  \chi_{S^-} \}$ generate 
the character group of   $\C^* \times \mu_a$. Moreover, we have 
\[  \chi_D^{ap} (t, \zeta) \chi_{S^+}^m (t, \zeta) = 
 \chi_D^{-aq} (t, \zeta) \chi_{S^-}^m (t, \zeta) = 1 \;\; \forall t \in \C^*, 
\; \forall \zeta\in \mu_a. \]
This implies the following two relations in ${\rm Cl}(E_{h,m})$:
\[ ap[{D}] + m [S^+] =  -aq[{D}] + m [S^-] = 0. \]
Consider two natural surjective homomorphisms
\[ \psi^+\; : \; \Z^2 \to {\rm Cl}(E_{h,m}), \; (k_1, k_2) \mapsto 
k_1[{D}] +  k_2 [S^+], \]
\[ \psi^-\; : \; \Z^2 \to {\rm Cl}(E_{h,m}), \; (k_1, k_2) \mapsto 
k_1[{D}] +  k_2 [S^-]. \]
Then 
$${\rm Ker}\, \psi^+ = \langle (ap, m) \rangle, \; \;
{\rm Ker}\, \psi^- = \langle (-aq, m) \rangle, $$
because each of two elements $(p, k), (-q,k) \in \Z^2$ generates 
a direct summand of $\Z^2$, and,   by 
$ka =m$, we have 
\[ \Z^2/\langle (pa,m) \rangle \cong \Z \oplus C_a \cong \Z^2/\langle 
(-qa,m) \rangle. \]
 \hfill $\Box$ 
\medskip

\section{SL(2)-equivariant flips} 

Let us start with toric $SL(2)$-equivariant  flips. 
It is known that if $m = a(q-p)$ 
then  the toric variety $E_{h,m}$ is isomorphic to the 
closure of the orbit of the highest vector in the 
irreducible $SL(2)\times SL(2)$-module  $V_{ap} \otimes V_{aq}$ 
\cite[Prop.2]{Pa92}.  In this case, 
$E_{h,m}$ is isomorphic to the 
affine cone in $V_{ap} \otimes V_{aq} \cong 
\C^{(ap+1)\times (aq+1)}$ with vertex $0$
over the projective 
embedding of $\P^1 \times \P^1$ into a projective space 
by the global sections of the 
ample sheaf ${\mathcal O}(ap, aq)$. The closure $D$ of the 
$2$-dimensional 
$SL(2)$-orbit ${\mathcal D}$ in  $E_{h,m}$ is isomorphic to the 
affine cone over $a(p+q)$-th Veronese embedding of 
$\P^1$ considered as diagonal in $\P^1 \times 
\P^1$. If $e_1,e_2,e_3$ is a standard basis of $\R^3$ then 
the toric variety  $E_{h,m}$ is defined by the cone 
$\sigma = \sum_{i =1}^4 \R_{\geq 0} v_i$ where 
$$v_1 = e_1, \; v_2 = -e_1 + aq e_3, \; 
v_3 = e_2, v_4 = -e_2 + ap e_3, $$ 
i.e., $v_1, v_2, v_3, v_4$ satisfy the equation 
$pv_1 + pv_2 = qv_3 + qv_4$.

Let ${E_{h,m}'}$ 
be the blow up of $0 \in E_{h,m} \subset \C^{(ap+1)\times (aq+1)}$. 
It corresponds to the subdivion of $\sigma$ into $4$ simplicial cones
having a new common ray $\R_{\geq 0} v_5$ $(v_5 = e_3)$ 
and generated by the following 4 sets of lattice vectors 
\[ \{v_1,v_3, v_5 \}, \;  \{v_2,v_3, v_5 \}, \; 
 \{v_2,v_4, v_5 \},\;  \{v_4, v_1, v_5 \}. \] 
The exceptional divisor $D'$
over $0$ corresponding  to the new lattice vector 
$v_5$ is isomorphic to  $\P^1 \times \P^1$. Moreover, the whole variety 
 ${E_{h,m}'}$ is smooth and 
can be considered 
as a line bundle of bidegree $(-aq,-ap)$ over $\P^1 \times \P^1$. 
Consider two $2$-dimensional simplical cones 
\[ \sigma^+ =  \R_{\geq 0} v_3 +   \R_{\geq 0} v_4, \; \mbox{\rm and} \; 
 \sigma^- =  \R_{\geq 0} v_1 +   \R_{\geq 0} v_2. \]
There exist two different subdivisons of $\sigma$ into pairs of simplicial 
cones 
\[ \sigma = (\R_{\geq 0} v_1 +   \sigma^+) \cup 
 (\R_{\geq 0} v_2 +  \sigma^+) \; \mbox{\rm and} \; 
\sigma = (\R_{\geq 0} v_3 +   \sigma^- ) \cup 
 (\R_{\geq 0} v_4 +  \sigma^-). \]
We denote toric varieties corresponding two these subdivisions by 
$ E_{h,m}^-$ and $ E_{h,m}^+$ respectively. 
Then one obtains the following diagram of toric morphisms:  
\[
\xymatrix{ & {E_{h,m}'}  \ar[dr]^{\gamma^+}   \ar[dl]_{\gamma^-} &\\
 E_{h,m}^-  \ar[dr]_{\varphi^-}& & E_{h,m}^+  \ar[dl]^{\varphi^+}\\
&E_{h,m}  & }
\]
The morphisms $\gamma^-$ and $\gamma^+$ restriced to 
$D'$ are projections of $\P^1 \times \P^1$ onto first and second factors.  
We denote by $C^-$ (reps. $C^+$) the $\gamma^-$-image (resp. $\gamma^+$-image) 
of $D'$ in $  E_{h,m}^-$ (resp.  $E_{h,m}^+$). 
Then singularities along  $C^-$ (reps. along $C^+$) are determined by 
the $2$-dimensional cone $\sigma^-$ (resp.  $\sigma^+$). 
The relations
\[ v_3 + v_4 = apv_5, \;\;   v_1 + v_2 = aqv_5 \]
show that the $2$-dimensional affine toric variety 
$X_{\sigma^-}$ (resp. $X_{\sigma^+}$)  is an  affine cone over 
$\P^1$ embedded by ${\mathcal O}(ap)$ (resp.  by ${\mathcal O}(aq)$) 
to $\P^{ap}$ (resp.  $\P^{aq}$). 
By  $1 \leq p <q$, we obtain that   $E_{h,m}^-$ is always singular and 
$E_{h,m}^+$ is nonsingular if and only if $ap=1$. Simple calculations in 
Chow rings of toric varieties  $E_{h,m}^-$ and  $E_{h,m}^+$ show that  
\[ C^- \cdot  K_{E_{h,m}^-} = 
\frac{2(p-q)}{aq^2}  < 0, \;\; \;\;  
C^+ \cdot  K_{E_{h,m}^+} = \frac{2(q-p)}{ap^2} > 0. \]
So the birational map 
$$\xymatrix{ E_{h,m}^- \ar@{-->}[r] &  E_{h,m}^+ }$$ is a  
toric flip.     
\medskip

Now we consider a general case for an affine $SL(2)$-variety  $E_{h,m}$. 
Let us begin with the calculation of the canonical class of an arbitrary  
$SL(2)$-variety $E_{h,m}$ which has been done  
by Panyushev in \cite[Prop.4 and 5]{Pa92}: 

\begin{proposition} 
For any normal affine $SL(2)$-variety $E_{h,m}$, one has 
\[ K_{E_{h,m}} = -(1 + b) [D]. \]
\label{can-class}
\end{proposition}

\noindent
{\em Proof.} Using the description of $E_{h,m}$ as a categorical 
quotient $H_b/\!\!/G$ of the hypersurface $H_b \subset \C^5$, we 
can consider $E_{h,m}$ as a hypersurface in the $4$-dimensional 
affine toric variety ${\mathcal T}_{h,m}:= 
\C^5/ \!\!/G$. It is well-known 
that the canonical divisor 
of any toric variety consists of irreducible divisors in the complement 
to the open torus orbit taken with the multiplicity $-1$. If we 
consider $Y_0, X_1, X_2, X_3, X_4$ as homogeneous coordinates 
of the toric variety   ${\mathcal T}_{h,m}$, then the  canonical class 
of  ${\mathcal T}_{h,m}$ corresponds to the character $\chi\, : \, G \to \C^*$ 
\[  \chi(t, \zeta) = t^{-k} (t^{p}\zeta)^2 (t^{-q}\zeta^{-1})^2 = 
t^{-k +2p-2q}. \]
On the other hand, $G$ acts on the polynomial $Y_0^b - X_1X_4 + X_2X_3$ 
by the character 
\[  \chi'(t, \zeta) = t^{q-p}. \]
Therefore, by adjunction formula, the canonical class of 
 $E_{h,m}$ corresponds to the character $\chi^+ = \chi + \chi'$: 
\[ \chi^+(t, \zeta)  = t^{-k + p -q}. \]
Since the class $[D] \in {\rm Cl}(E_{h,m})$ 
is defined by  the character 
$\chi_D(t, \zeta) = t^k$, we obtain that 
\[   K_{E_{h,m}} = \frac{ -k + p -q}{k} [D] = - (1+b)[D]. \]
\hfill $\Box$

\begin{proposition} 
Let $L^+$ be the trivial line bundle over $H_b$ together with the 
linearization corresponding to the character $\chi^+$, then 
\[  H_b^{ss}(L^+) = U^+ = H_b \setminus \{ X_1 = X_2 = 0\}. \]
\end{proposition} 

\noindent
{\em Proof.} The space $\Gamma(H_b, (L^+)^{\otimes n})^G$ consists 
of all regular functions $f$ on $H_b$ such that 
$f(gx) = (\chi^+(g))^n f(x)$. It is easy to see that   
$\Gamma(H_b, (L^+)^{\otimes n})^G$ is generated as a $\C$-vector space 
by restrictions of monomials $Y_0^{k_0} X_1^{k_1}  X_2^{k_2}  X_3^{k_3} 
X_4^{k_4}$ satisfying the above  homogeneity condition, i.e. , 
\[  t^{k_0k -k_1p - k_2p + k_3 + q k_4q} \zeta^{-k_1 -k_2 + k_3 + k_4} = 
t^{n(-k + p-q)} \;\; \forall t \in \C^*, \;\; \forall \zeta \in \mu_a. \]
The last condition implies $a | ( k_3 + k_4 - k_1 -k_2)$ and 
\[ k_0k -k_1p - k_2p + k_3 q +  k_4q = n(-k + p-q). \]
Since  $n(-k + p-q) < 0$ and $k_i \geq 0$ $(0 \leq i \leq 4)$, 
we obtain that at least one of the integers $k_1$ and $k_2$ must 
be positive, i.e., all monomials 
$Y_0^{k_0} X_1^{k_1}  X_2^{k_2}  X_3^{k_3} 
X_4^{k_4} \in \Gamma(H_b, (L^+)^{\otimes n})^G$ 
vanish on the subset  $\{ X_1 = X_2 = 0\} \cap H_b$. On the other hand, 
if at least one of two coordinates $X_1$ and $X_2$ of a point $x \in H_b$ 
is not zero, then one of the monomials 
\[  X_1^{q-p +k} , \;\; X_2^{q-p+k} \in \Gamma(H_b, (L^+)^{\otimes p})^G \]
does not vanish in $x$.  Hence,  $H_b^{ss}(L^+) = U^+$.
\hfill $\Box$

\begin{proposition} 
Let $L^-$ be the trivial line bundle over $H_b$ together with the 
linearization corresponding to the character $\chi^- = - \chi^+$, then 
\[  H_b^{ss}(L^-) = U^- = H_b \setminus \{ X_3 = X_4 = 0\}. \]
\end{proposition} 

\noindent
{\em Proof.} The condition $f(gx) = (\chi^-(g))^n f(x)$ for  
a monomial  $$f = Y_0^{k_0} X_1^{k_1}  X_2^{k_2}  X_3^{k_3} 
X_4^{k_4} \in  \Gamma(H_b, (L^-)^{\otimes n})^G$$ implies that 
\[  t^{k_0k -k_1p - k_2p + k_3 q + k_4q} \zeta^{-k_1 -k_2 + k_3 + k_4} = 
t^{n(k + q-p)} \;\; \forall t \in \C^*, \;\; \forall \zeta \in \mu_a. \]
Since $n(k + q-p) >0$, we obtain that at least one of three integers 
$k_0, k_3, k_4$ must be  positive. Therefore, 
all monomials 
$Y_0^{k_0} X_1^{k_1}  X_2^{k_2}  X_3^{k_3} 
X_4^{k_4} \in \Gamma(H_b, (L^-)^{\otimes n})^G$ 
vanish on the subset  $\{ Y_0= X_3 = X_4 = 0\} \cap H_b = 
\{ X_3 = X_4 = 0\} \cap H_b$. 
 On the other hand, 
if at least one of two coordinates $X_3$ and $X_4$ of a point $x \in H_b$ 
is not zero, then one of the monomials 
\[  X_3^{q-p +k} , \;\; X_4^{q-p+k} \in \Gamma(H_b, (L^-)^{\otimes q})^G \]
does not vanish in $x$. Hence,  $H_b^{ss}(L^-) = U^-$. 

\hfill $\Box$

\begin{theorem} 
Define 
\[ E_{h,m}^- :=  H_b^{ss}(L^-)/\!\!/G , 
\;\;   E_{h,m}^+ :=  H_b^{ss}(L^+)/\!\!/G.\]
Then the  open embeddings 
\[  H_b^{ss}(L^-) = U^- \subset H_b, \;\;  H_b^{ss}(L^+) = U^+ \subset H_b, \]
define two natural birational morphisms 
\[ \varphi^-\; : \;  E_{h,m}^- \to  E_{h,m} , \;\;  \varphi^+\; : \;  
E_{h,m}^+ \to  E_{h,m}, \]
and the  $SL(2)$-equivariant flip 
\[
\xymatrix{ E_{h,m}^- \ar@{-->}[rr] \ar[dr]_{\varphi^-}& & E_{h,m}^+  
\ar[dl]^{\varphi^+}\\
&E_{h,m}  & }
\]
\label{flip1}
\end{theorem} 

\noindent
{\em Proof.} The statement follows immediately from the isomorphisms
\[  E_{h,m}^- \cong {\rm Proj} \, \bigoplus_{n \geq 0}  \Gamma(H_b, 
(L^-)^{\otimes n})^G \cong 
{\rm Proj} \, \bigoplus_{n \geq 0}  \Gamma(E_{h,m}, 
{\mathcal O}(-nK_{E_{h,m}}) \]
and 
\[  E_{h,m}^+ \cong {\rm Proj} \, \bigoplus_{n \geq 0}  \Gamma(H_b, 
(L^+)^{\otimes n})^G \cong 
{\rm Proj} \, \bigoplus_{n \geq 0}  \Gamma(E_{h,m}, 
{\mathcal O}(nK_{E_{h,m}}). \]

\hfill $\Box$

\begin{corollary}
One has the following isomorphisms: 
\[   E_{h,m}^- \cong {\rm Proj} \, \bigoplus_{n \geq 0}  \Gamma(E_{h,m}, 
{\mathcal O}(-nD)),  \;\; \; 
 E_{h,m}^+ \cong {\rm Proj} \, \bigoplus_{n \geq 0}  \Gamma(E_{h,m}, 
{\mathcal O}(nD)). \]
\end{corollary}

\noindent
{\em Proof.} These isomorphisms follow from 
the equation $K_{ E_{h,m}} = -(1+b)[D]$ (\ref{can-class}) 
and from the isomorphism 
\[ {\rm Proj}\, \bigoplus_{n \geq 0} R_n \cong   
{\rm Proj}\, \bigoplus_{n \geq 0} R_{nl} \]
for any notherian graded ring   $R= \bigoplus_{n \geq 0} R_n$ and for 
any positive integer $l$. 
\hfill $\Box$ 
\bigskip

In order to describe the geometry  $E_{h,m}^-$ and $ E_{h,m}^+$ in more detail 
we need two $2$-dimensional affine varieties $S^+$ and 
$S^-$ having regular $B$-actions (see also \ref{cl-div}).

\begin{proposition} Let   $S^+ \subset E_{h,m}$ be the closure 
of an $B$-orbit obtained as categorical quotient of $W^+ := 
 H_{q-p} \cap \{ X_2 =0 \}$ by $G_0 \times G_m$. Then  $S^+$ 
is isomorphic to the normal affine toric surface 
${\rm Spec}\, \C[ M_{h,m}^+ ]$.  
\label{X2}
\end{proposition} 

\noindent
{\em Proof.}  We note that  $W^+= H_{q-p} \cap \{ X_2 =0 \} \subset \C^5$ 
is a $3$-dimesional affine 
toric variety which is a product of $\C$ and a $2$-dimensional 
affine toric variety defined by the binomial 
equation $X_0^{q-p} = X_1 X_4$. Let us 
compute the categorical quotient $W^+/\!\!/ G_0$. Since $G_0$ acts on 
$X_0, X_1, X_3, X_4$ by $diag (t, t^{-p}, t^{q}, t^q)$ for every nonconstant 
$G_0$-invariant monomial 
$X_0^{k_0} X_1^{k_1}  X_3^{k_3}  X_4^{k_4}$ $(k_i \in \Z_{\geq 0})$ 
the condition $k_0 - pk_1 + qk_3 + qk_4 =0$ implies $k_1 >0$. 
If at the same time $k_4 > 0$, then 
\[ X_0^{k_0} X_1^{k_1}  X_3^{k_3}  X_4^{k_4} -  
X_0^{k_0+ q-p} X_1^{k_1-1}  X_3^{k_3}  X_4^{k_4-1} \in I(W^-). \]
Using  the equation $X_0^{q-p} = X_1 X_4$ several times, 
we can get another monomial
$X_0^{k_0'} X_1^{k_1'}  X_3^{k_3'}$ such that   
$X_0^{k_0} X_1^{k_1}  X_3^{k_3}  X_4^{k_4} - X_0^{k_0'} X_1^{k_1'}  
X_3^{k_3'} \in I(W^+)$, i.e.,  
vanish on $W^+$. Therefore, the 
coordinate ring of $W^+/\!\! /G_0$ contains a 
$\C$-basis consisting of all
 $G_0$-invariant 
monomials in $X_0, X_1, X_3$. 
These monomials have form  $X_0^{pk_1 -qk_3} X_1^{k_1} X_3^{k_3} = 
X^{k_1} Y^{k_3}$ 
where $pk_1-qk_3 \geq 0$ (i.e. $(k_1,k_3) \in M_{h,1}^+)$). So the coordinate 
ring of $S^+ = W^+ / \!\!/(G_0 \times G_m)$ has a $\C$-basis consisting 
of $G_m$-invariants monomials $X^{k_1} Y^{k_3} = 
X_0^{pk_1 -qk_2} X_1^{k_1} X_3^{k_3}$ which correspond 
to lattice points 
$(k_1,k_3) \in M_{h,m}^+ = M_{h,1}^+ \cap \{ (k_1,k_3) \in \Z^2_{\geq 0}\,  
:\, m | (k_1-k_3) 
\}$, i.e.  $S^+ \cong  {\rm Spec}\, \C[ M_{h,m}^+ ]$.     \hfill $\Box$

\begin{proposition} Let   $S^- \subset E_{h,m}$ be the closure 
of an $B$-orbit 
 obtained as categorical quotient of $W^- := 
 H_{q-p} \cap \{ X_4 =0 \}$ by $G_0 \times G_m$.
Then  $S^-$ 
is isomorphic to the normal affine toric surface  ${\rm Spec}\, 
\C[ M_{h,m}^- ]$, 
where the monoid $ M_{h,m}^- \subset \Z^2$ (see Figure 2) 
is defined as follows: 
\[ M_{h,m}^-  := \{ (i,j) \in \Z^2 \; : \; j \leq h i, \; 
i \geq 0, \; m | (i-j) \}.  \]
\label{X4}
\end{proposition} 
\begin{figure}[htbp]\label{bild2}
\begin{center}
\includegraphics[width=0.7\textwidth]{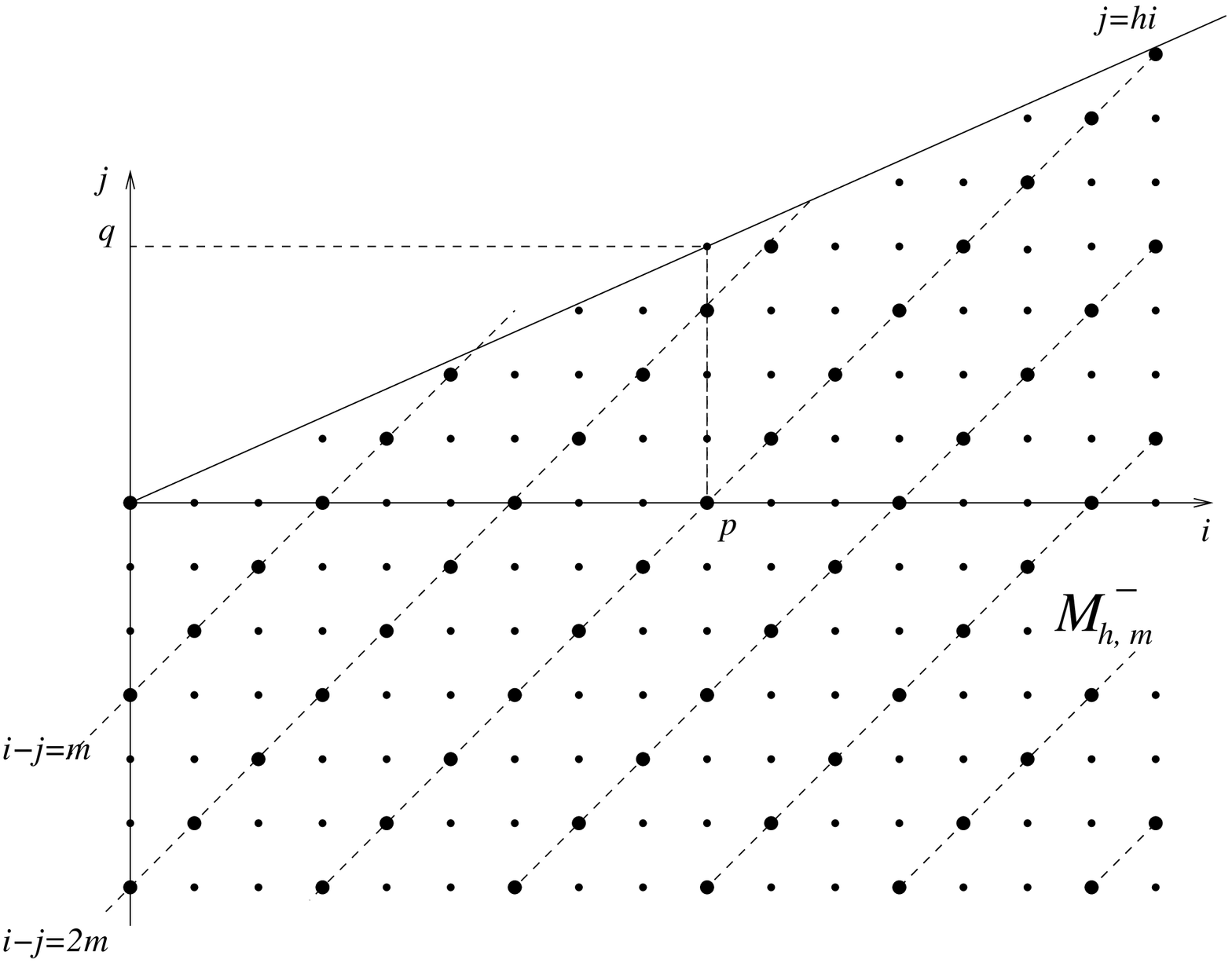}
\end{center}
\caption{}
\end{figure}
\noindent
{\em Proof.}  We note that  $W^-= H_{q-p} \cap \{ X_4 =0 \} \subset \C^5$ 
is a $3$-dimesional 
toric variety which is a product of $\C$ and a $2$-dimensional 
toric variety defined by the equation $X_0^{q-p} = -X_2 X_3$. Again 
the computation  of the categorical quotient 
$W^- /\!\!/ G_0$ reduces to finding all $G_0$-invariant
monomials $X_0^{k_0} X_1^{k_1}  X_2^{k_2}  X_3^{k_3}$.  
Under the condition 
$X_0^{q-p} = -X_2 X_3$ we can assume that at least one   
of two variables $X_2$, or $X_3$ does not appear in 
$X_0^{k_0} X_1^{k_1}  X_2^{k_2}  X_3^{k_3}$ 
(i.e., $k_2 =0$ ,or $k_3 = 0$). If $k_2 =0$, then we come to the same 
situation as in \ref{X2} and obtain $G_0$-invariant 
monomials   $X^{k_1} Y^{k_3} = X_0^{pk_1 -qk_3} 
X_1^{k_1} X_3^{k_3}$ $(k_1,k_3) \in M_{h,1}^+$. 
If $k_3 =0$, then we obtain $G_0$-invariant monomials  
$X_0^{pk_1 +pk_2} X_1^{k_1} X_2^{k_2} = 
X^{k_1} Z^{k_2}$,  $(k_1,k_2 \in \Z_{\geq 0})$. 
The equation $X_0^{q-p} = -X_2 X_3$ implies  that on $W^-/\!\!/G_0$ 
we have $YZ= X_0^{-q}X_2 X_0^p X_3 = -1$.  So in case $k_2 =0$ 
we obtain the monomials in $X^{k_1} (Y^{-1})^{k_2}$,  $(k_1,k_2 
\in \Z_{\geq 0})$. 
Unifying both cases, we get all $G_0$-invariant 
monomials $X^iY^j$,  $(i,j ) \in \widetilde{M}_h^1$. The action of the
finite group $G_m$ on $X$ and $Y$ gives rise to 
an additional restiriction: $m | (i-j)$. Therefore, 
$G_0 \times G_m$-invarinant monomials can be identified with 
the set of all lattice points $(i,j) \in M_{h,m}^-$. 

\hfill $\Box$

\begin{remark} 
{\rm If $m = a(q-p)$ (i.e. $E_{h,m}$ is toric),  then 
$ S^- \cong X_{\sigma^-}$ and  $S^+ \cong X_{\sigma^+}$, where
$\sigma^-$ and $\sigma^+$ are $2$-dimensional cones as above. 
} 
\end{remark} 

\begin{definition} 
{\rm Let $S$ be an algebraic surface with a regular action $B \times S \to S$ 
of a Borel subgroup $B \subset SL(2)$. We denote by 
$SL(2) \times_B S$
the $SL(2)$-variety $(SL(2) \times S)/B$, where $B$ is considered to act on 
$SL(2)$ by right multiplication: 
\[ \left( \begin{pmatrix} X & Y \\ Z & W \end{pmatrix}, 
\begin{pmatrix} t & u \\ 0 & t^{-1} \end{pmatrix} \right) \mapsto 
\begin{pmatrix} X & Y \\ Z & W \end{pmatrix}\cdot 
 \begin{pmatrix}  t & u \\ 0 & t^{-1} \end{pmatrix}^{-1}, \; \;\; 
 \;\; \begin{pmatrix} X & Y \\ Z & W \end{pmatrix} \in 
SL(2).\]  
} 
\end{definition} 

\begin{theorem} One has the following isomorphisms
\[ E_{h,m}^- \cong   SL(2) \times_B S^- , 
\;\;   E_{h,m}^+ \cong  SL(2) \times_B S^+.\]
\end{theorem}

\noindent
{\em Proof.}  Since $U^+ = H_b \setminus \{ X_1 = X_2 = 0 \}$ and $G$ acts 
on $(X_1, X_2)$ by scalar matrices,  
we obtain a natural $SL(2)$-equivariant morphism 
$$ \alpha^+\; : \;  E_{h,m}^+ \cong   U^+/\!\!/ G 
\to \P^1, \; \; (Y_0, X_1, X_2, X_3, X_4) \mapsto 
(X_1: X_2) $$
Analogously, we obtain a natural $SL(2)$-equivariant morphism
$$ \alpha^-\; : \;  E_{h,m}^- \cong  U^-/\!\!/ G 
\to \P^1, \; \; (Y_0, X_1, X_2, X_3, X_4) \mapsto 
(X_3: X_4). $$
By \ref{X2} and \ref{X4}, we have 
$$S^+ = (\alpha^+)^{-1}( 1:0), \;\;  S^- = (\alpha^+)^{-1}( 1:0). $$
Since the morphisms  $\alpha^+$ and $\alpha^-$ are $SL(2)$-equivariant 
and $SL(2)$ acts transitively on $\P^1$, we have 
\[ S^+ \cong  (\alpha^+)^{-1}(z), \;  S^- \cong  (\alpha^-)^{-1}(z) 
\; \; \forall z \in \P^1, \]
i.e., $E^+_{h,m}$ (resp. $E_{h,m}^-$) is a fibration over $\P^1$ 
with fiber $S^+$ (resp. $S^-$).  
On the other hand, the projection $SL(2) \times S^{\pm} \to SL(2)$
defines two natural morphisms $SL(2)$-equivariant morphisms  
\[ \pi^+ \; : \;  SL(2) \times_B S^+ \to SL(2)/B \cong \P^1, \]
\[ \pi^- \; : \;  SL(2) \times_B S^- \to SL(2)/B \cong \P^1, \]
such that $ SL(2) \times_B S^+$ (resp.  $SL(2) \times_B S^-$) 
is a fibration over $\P^1$ with 
fiber $S^+$ (resp. $S^-$). Consider  the 
morphisms 
\[  \tilde{\beta}^+ \; : \; SL(2) \times S^+  \to U^+/\!\!/G, \;\;  
 \tilde{\beta}^- \; : \; SL(2) \times S^-  \to U^-/\!\!/G  \]
defined by 
\[   \tilde{\beta}^\pm ( g, x) = gx, \;\; \forall g \in SL(2), \; 
\forall x \in S^{\pm} = (\alpha^{\pm})^{-1}(1:0). \]   
Since $ \tilde{\beta}^\pm ( gb^{-1}, bx) = \tilde{\beta}^\pm ( g, x) = 
gx$ $\forall b \in B$, the morphism  
$\tilde{\beta}^\pm$ descends
to a morphism 
\[  {\beta}^{\pm} \; : \; SL(2) \times_B S^{\pm} \to  
U^{\pm}/\!\!/G \cong E_{h,m}^{\pm}. \] 
The latter  is an isomorphism, because  ${\beta}^{\pm}$ is 
$SL(2)$-equivariant and it maps isomorphically  
the fiber of $\pi^{\pm}$ over $(1:0)$ to the fiber of $\alpha^{\pm}$ over 
the $B$-fixed point $[B] \in SL(2)/B$. 
\hfill $\Box$

\begin{remark} 
{\rm Since  the monoid $M_{h,m}^+ $ is submonoid of the monoid 
$M_{h,m}^-$ we obtain a birational morphism $\psi\; : \; 
S^- \to S^+$ 
of $2$-dimensional normal affine toric varieties  
$S^-$ and $S^+$. However, $\psi$ is not $B$-equivariant, because
an  element 
$$\begin{pmatrix} t & u \\ 0 & t^{-1} \end{pmatrix} \in B $$
sends $X^i Y^j \in \C[ M_{h}^m ]$ to $$(tX)^i (tY + uX^{-1})^j \in  
\C[ M_{h}^m ]$$ 
and sends   $X^r Y^s \in \C[M_{h,m}^- ]$  to 
$$( tX - uY^{-1} )^r (tY)^s \in  \C[M_{h,m}^- ].$$ 
This is the reason why 
there is no any birational $SL(2)$-equivariant
morphism from  $SL(2) \times_B S^-$  to   $SL(2) \times_B S^+$, 
but only a flip.} 
\end{remark}

\begin{remark} 
{\rm Let $E_{h,m} \hookrightarrow V$ be  a closed embedding, where 
$V$ is an affine space isomorphic to $V_{i_1 +j_1} \oplus \cdots \oplus 
V_{i_r +j_r}$ (see \ref{h<1}). We define a $\C^*$-action on $V$ such that 
$t \in \C^*$ acts by multiplication with $t^{j - i}$ on  $V_{i +j}$. 
Since this $\C^*$-action commutes with the $SL(2)$-action, the affine 
variety $E_{h,m} \subset V$ remains  invariant under the $\C^*$-action. 
Consider the weighted blow up $\delta\, : \, 
\widetilde{V} \to V$ of $0 \in V$ with respect to weights of 
this  $\C^*$-action. 
The birational pullback 
of  $E_{h,m}$ under $\delta\; : \; 
\widetilde{V} \to V$ is a  
$SL(2)$-variety ${E_{h,m}'}$ together 
with surjective morphisms  
$\gamma^{-}\, : \, {E_{h,m}'} \to 
{E_{h,m}^-}$ and  $\gamma^{+}\, : \, {E_{h,m}'} \to E_{h,m}^+$ 
such that the following 
diagram commutes
\[ \xymatrix{ & {E_{h,m}'}  \ar[dr]^{\gamma^+}   \ar[dl]_{\gamma^-} &\\
 E_{h,m}^-  \ar@{-->}[rr] & & E_{h,m}^+ }  
 \]
The variety  ${E_{h,m}'}$ contains two  $SL(2)$-invariant 
divisors $D' \cong  \P^1 \times \P^1$ and $\widetilde{D} := \delta^*(D)$ 
whose  intersection  $C= D' \cap \widetilde{D} \cong 
\P^1$ is the unique $1$-dimensional closed $SL(2)$-orbit in  ${E_{h,m}'}$.  
The morphism $\gamma^{\pm}$ contracts  $D'$ to $C^{\pm} \subset E_{h,m}^{\pm}$.
The divisor $D'$ corresponds to the $SL(2)$-invariant discrete 
valuation of the function field $\C( SL(2))$ 
defined by above $\C^*$-action on  $E_{h,m}$ such that $\C(D')$ is 
the  $\C^*$-invariant 
subfield $\C( SL(2))^{\C^*} \cong \C( SL(2)/\C^*)$.    
We note that  the $SL(2)$-variety ${E_{h,m}'}$ has also 
a toroidal structure, i.e., 
along the closed $1$-dimensional $SL(2)$-orbit $C$,  
it is locally isomorphic to 
a product of an affine line $\A^1$ and a $2$-dimensional 
affine toric surface $S'$  which  is isomorphic to 
${\rm Spec}\, \C[M_{h,m}']$ where  
\[ M_{h,m}':= \{ (i,j) \in \Z^2\; : \; pj - qi \geq 0, j - i \in  
m \Z_{\geq 0} \}. \]
In particular, $S'\cong 
\A^2/\mu_b$ and ${E_{h,m}'}$ is nonsingular along $C$ 
if and only if $b =1$, i.e., iff $E_{h,m}$ is toric.
} 
\label{blow-up}
\end{remark}

\begin{proposition} 
The canonical divisor of $E_{h,m}^\pm $ has the following 
intersection numbers with the $1$-dimensional $SL(2)$-orbits 
$C^{\pm} \subset  E_{h,m}^\pm$:
\[  K_{E_{h,m}^-} \cdot  C^-  
= - \frac{(1+b)k}{aq^2},\;\;  
K_{E_{h,m}^+} \cdot  C^+  = \frac{(1+b)k}{ap^2}. \]   
\end{proposition} 

\noindent
{\em Proof.} 
Since $E_{h,m}, E_{h,m}^-$ and $E_{h,m}^+$ 
have the same divisor class group, we can use 
\ref{cl-div} and obtain that 
\[ ap[D] + m[S^+] = 0 \in {\rm Cl}(E_{h,m}^+). \]
The divisor  $S^+ \subset E_{h,m}^+$ intersects the curve $C^+$ 
transversally, but this intersection point is an isolated 
cyclic quotient singularity of type $A_{ap-1}$ in $S^+$. Therefore, we 
have 
$S^+ \cdot C^+ = \frac{1}{ap}$ and   
$$D \cdot C^+  = - 
\left(\frac{m}{ap} S^+ \right) \cdot C^+ = - \frac{k}{ap^2}. $$
By \ref{can-class}, we get 
\[ K_{E_{h,m}^+} \cdot  C^+  = \frac{(1+b)k}{ap^2}.\]
Similarly, the intersection point of $C^-$ and  $S^- \subset E_{h,m}^-$
is an isolated cyclic quotient singularity of type $A_{aq-1}$ in $S^-$.
 Therefore, we 
have $S^- \cdot C^- = \frac{1}{aq}$ and,  by $$-aq[D] + m[S^-] = 0 
\in {\rm Cl}(E_{h,m}^-),$$ we obtain   
$$D \cdot C^-  =  \left(
\frac{m}{aq}  S^- \right) \cdot C^- =  \frac{k}{aq^2}. $$
By \ref{can-class}, this implies 
\[ K_{E_{h,m}^-} \cdot  C^-  = - \frac{(1+b)k}{aq^2}.\]
 \hfill $\Box$

\begin{remark} 
{\rm In  \cite[Section 9]{LV83}  Luna and Vust gave a description of 
an  arbitrary quasihomogeneous normal $SL(2)$-embedding by a  special 
combinatorial diagram (so called ``marked 
hedgehog''). If we apply this language to the discription 
of  quasiprojective varieties  $E_{h,1}^-$ and   $E_{h,1}^+$ 
(here we assume $m =1$), then we find  a  difference between the  
using of signs $+$ and 
$-$ in our paper and in \cite[Section 9]{LV83}. For instance,  
the $1$-dimensional orbit  $C^+ \subset E_{h,1}^+$ 
is considered  in \cite{LV83} as an $SL(2)$-orbit of type  
 $l_{-}(D,h)$. Similarly,  
the $1$-dimensional orbit $C^- \subset E_{h,1}^-$ 
is considered in  in \cite{LV83} 
as orbit of type  $l_{+}(D,h)$. 
} 
\end{remark}

\bigskip

\section{Flips and spherical varieties}

Let us consider the $\C^*$-action on $H_b$ defined by the 
diagonal matrices $$diag(1, s^{-1}, s^{-1}, s, s), \;\; s \in \C^*.$$ 
We note that this   $\C^*$-action commutes with the 
$SL(2)$-action and with the action of 
$G= G_0' \times G_a$. So we obtain a natural $\C^*$-action 
on the categorical quotient $H_b / \!\!/ G \cong E_{h,m}$ which commutes with 
the $SL(2)$-action. We note that 
this $\C^*$-action has been already constructed in \ref{blow-up} using 
a closed embedding $E_{h,m} \hookrightarrow V$. 
This allows to consider $E_{h,m}$ as an affine  
$SL(2) \times \C^*$-variety. 

\begin{proposition} 
The affine variety  $E_{h,m}$ is spherical with respect 
to the above  $SL(2) \times \C^*$-action. 
\end{proposition} 

\noindent
{\em Proof.} The open subset ${\mathcal U}= 
(H_b \cap  \{ Y_0 \neq 0\})/G \subset E_{h,m}$ is obviously  
$SL(2) \times \C^*$-invariant. Since  $SL(2) \times \C^*$ acts trasitively on 
 ${\mathcal U}$, we have  ${\mathcal U} \cong (SL(2) \times \C^*)/H$ 
for some closed subgroup $H \subset SL(2) \times \C^*$. It is easy 
to see that  $$(H_b \cap  \{ Y_0 X_2 X_4 \neq 0\})/G \subset {\mathcal U}$$ 
is an open dense orbit of the $3$-dimensional Borel subgroup $\widetilde{B}:= 
B \times \C^*$ in $SL(2) 
\times \C^*$. Hence,   $E_{h,m}$ is a spherical embedding corresponding to 
the spherical homogeneous space $(SL(2) \times \C^*)/H$. 
\hfill $\Box$    
\bigskip 

\begin{remark} 
{\rm There exists one more way to define the same $\C^*$ on 
$E_{h,m}$. We identify $\C^*$ with the maximal torus $T \subset SL(2)$ 
which acts on $SL(2)$ by right multiplication. Then this action extends 
to a regular action on $E_{h,m}$ and commutes with the $SL(2)$-action 
by left multiplication so that we obtain a regular action of 
$SL(2) \times T$ on $E_{h,m}$.  By \cite[]{Kr84}, even a more general 
statement is true: $E_{h,m}$ admits a regular action of $SL(2) \times B$.     
} 
\end{remark} 

If we identify ${\mathcal U}$ with $SL(2)/C_m$ and consider  the subgroup 
$H \subset SL(2) \times \C^*$ as 
a stabilizer of the class of unit matrix in  $SL(2)/C_m$ then 
\[ H= \{ (diag(t,t^{-1}), t^m) \; : \; t \in \C^* \} \subset 
 SL(2) \times \C^*. \]
The lattice $\Lambda$ of rational $\widetilde{B}$-eigenfunctions  on  
${\mathcal U}$ (up to multiplication 
with a nonzero constant) consists of all Laurent monomials 
$Z^{i} W^{j} \in \C[SL(2)]^{\mu_m}$ 
such that $m | (i-j)$. Therefore,  $E_{h,m}$ is a spherical 
embedding of rank $2$. This rank equals also the minimal codimension 
of $U$-orbits in $E_{h,m}$ (we identify $U$ with the maximal 
unipotent subgroup in $SL(2) \times \C^*$).

In order to describe spherical varieties  $E_{h,m}$,  $E_{h,m}^+$, and 
$E_{h,m}^-$ by combinatorial data, we remark that 
they contain  exactly three $\widetilde{B}$-invariant divisors: 
\[ D = H_b \cap \{Y_0 =0 \}/\!\!/G, \; \;  S^+ = 
H_b \cap \{X_2 =0 \}/\!\!/G,  \; \;  S^- = 
H_b \cap \{X_4 =0 \}/\!\!/G. \]
The restrictions of the corresponding descrete valuations $
\C({\mathcal U})^* \to \Z$   
to the lattice  $\Lambda$ 
define lattice vectors $\rho, \rho^+, \rho^- \in \Lambda^*$ 
in  the dual space ${\mathcal Q}:= {\rm Hom}(\Lambda, \Q)$. 
We can consider $\rho^+, \rho^-$ as a $\Q$-basis of  ${\mathcal Q}$.  
Then the set of all $SL(2)\times \C^*$-invariant valuations generate 
so called {\em valuation cone} ${\mathcal V} \subset {\mathcal Q}$,   
${\mathcal V} = \{ x \rho^+ + y \rho^- \in {\mathcal Q}\; : \; 
x + y \leq 0 \}:$

\begin{center}
\begin{picture}(100,200)(10,20)
\put(50,160){\vector(1,0){40}}
\put(50,160){\vector(0,1){40}}
\put(50,160){\vector(1,-1){40}}
\put(50,160){\vector(2,-3){80}}
\put(135,45){$\rho$}
\put(95,165){$\rho^+$}
\put(55,204){$\rho^-$}
\put(95,125){$\rho'$}
\put(0,0){\thicklines  \put(50,160){\line(-1,1){50}}}
\put(50,160){\line(1,-1){100}}

\multiput(50,160)(5,-5){20}
{\line(-1,-1){30}}

\multiput(50,160)(-5,5){10}
{\line(-1,-1){30}}

\put(0,160){${\mathcal V}$}

\end{picture}

\end{center}

\noindent
The equations $Z= X_0^pX_2$, $W = X_0^{-q}X_4$ imply   
\[ \rho = p \rho^+ - q \rho^- \in  {\mathcal V}. \]
It is easy to see that $E_{h,m}$, $E_{h,m}^-$, and $E_{h,m}^+$ 
are simple spherical embeddings (i.e., they contain exactly one closed 
$SL(2) \times \C^*$-orbit of dimension $1$, or $0$). 
Therefore, they can be described by so called 
{\em colored cones} $({\mathcal C}, {\mathcal F})$, where ${\mathcal F}$ 
is a subset of $\{ \rho^+, \rho^- \}$ and 
${\mathcal C} \subset {\mathcal Q}$ is a strictly convex 
cone generated by ${\mathcal F}$ and $\rho$.  
More precisely we have: 
\[ {\mathcal C}(E_{h,m}) = \Q_{\geq 0} \rho +  \Q_{\geq 0} \rho^-, \;\; 
  {\mathcal F}(E_{h,m}) = \{ \rho^+, \rho^- \}, \]
\[  {\mathcal C}(E_{h,m}^-) = \Q_{\geq 0} \rho +  \Q_{\geq 0} \rho^+, \;\; 
  {\mathcal F}(E_{h,m}^-) = \{ \rho^+\}, \]
\[ {\mathcal C}(E_{h,m}^+) = \Q_{\geq 0} \rho +  \Q_{\geq 0} \rho^-, \;\; 
  {\mathcal F}(E_{h,m}^+) = \{ \rho^- \}, \]
Moreover, the spherical variety ${E_{h,m}'}$ is also simple. However,  
 ${E_{h,m}'}$  contains 
one more  $SL(2) \times \C^*$-invariant divisor $D'$  
such that the restrictions of 
the corresponding discrete valuations to $\Lambda$ defines a lattice 
vector   $\rho' =  \rho^+ - \rho^- \in  {\mathcal V}$. In this case, we have 
\[ {\mathcal C}({E_{h,m}'}) = 
\Q_{\geq 0} \rho +  \Q_{\geq 0} \rho', \;\; 
  {\mathcal F}({E_{h,m}'}) = 
\emptyset, \]

\begin{remark}
{\rm We note that birational morphisms $f\, : \, W \to W'$ of simple 
spherical varieties $W, W'$ where 
$f \in \{ \varphi^-, \varphi^+, \gamma^-, \gamma^+\}$ 
has an interpretation in terms of colors. In our situation, we see 
that the set of colors ${\mathcal F}(W')$ is strickly larger 
than  ${\mathcal F}(W')$.  

In particular, the birational morphism 
$\varphi^-\; : \; E_{h,m}^- \to  E_{h,m}$ combinatorially means 
that the cone ${\mathcal C}(E_{h,m}^-) = {\mathcal C}(E_{h,m})$ remains 
unchanged, but it gets an additional color $\rho^+$: 
 ${\mathcal F}(E_{h,m}) = {\mathcal F}(E_{h,m}^-) \cup  \{\rho^+\}$. 
On the other hand, the birational morphism 
 $\varphi^+\; : \; E_{h,m}^+ \to  E_{h,m}$ also adds  
 an additional color $\rho^-$:
 ${\mathcal F}(E_{h,m}) = {\mathcal F}(E_{h,m}^+) \cup  \{\rho^-\}$ such that 
 the color  $\rho^+$ becomes an interior point of   ${\mathcal C}(E_{h,m})$. 
This agree with a general description of Mori contractions 
in \cite[3.4, 4.4]{B93}.  
} 
\end{remark}

\begin{remark} 
{\rm According to Alexeev and Brion \cite{AB04}, every spherical 
$G$-variety ${\mathcal X}$ 
admits a flat degeneration to a toric variety ${\mathcal X_0}$.  
In general case, 
there exist several degenerations depending on different reduced 
decompositions of the longest element $w_0$ in the Weyl group of the 
reductive group $G$. However, in the case $G = SL(2) \times \C^*$ the 
choice of such a decomposition is unique.   
A simplest example of such a toric degeneration appears in the 
case ${\mathcal X} := SL(2)$ 
considered as a spherical homogeneous 
space of $SL(2) \times \C^*$. Then  ${\mathcal X}_0 =  
\{X_1X_4 - X_2X_3 =0\}$ is a singular affine $3$-dimensional 
toric quadric. The corresponding deformation  is 
${\mathcal X}_0 = \lim_{t \to 0} 
{\mathcal X}_t$ where $ {\mathcal X}_t:= \{X_1X_4 - X_2X_3 =t\}$. 
}
\label{tor-deg} 
\end{remark}

Let $T_{h,m}$ be the toric degeneration of $E_{h,m}$. Then  
\[ T_{h,m} := {\rm Spec}\, \C[ \widetilde{M}_{h,m}] \] 
where the semigroup 
\[   \widetilde{M}_{h,m} := \{ (i,j,k) \in \Z^3_{\geq 0} \; : \; 
m | (j-i), \; jp - qi \geq 0, \; i + j \geq k\}. \]
has surjective homomorphism $\pi\, : \, 
(i,j,k) \mapsto (i,j)$ onto 
$M_{h,m}^+$ where elements $(i,j)$ can be identified with the 
hightest vector $X^iY^j \in V_{i+j}$ and the lattice points 
$\pi^{-1}(i,j) \subset  \widetilde{M}_{h,m}$ correspond to 
the standard basis of $V_{i+j}$.  
So the toric degeneration  $T_{h,m}$ of $E_{h,m}$ is defined by a 
$3$-dimensional cone 
$$\sigma_0 = 
\R_{\geq 0} v_1 + \R_{\geq 0} v_2 + \R_{\geq 0} v_3 +  \R_{\geq 0} v_4$$
where $v_1 = (0,0,1)$, $v_2 = (1,1,-1)$, $v_3 = (0,1,0)$, $v_4 = (p,-q,0)$ 
satisfying the relation $$pv_1 + pv_2 = (p+q)v_3 + v_4.$$ 
In the notations of \cite{AB04}, the dual $3$-dimensional cone 
$\check{\sigma}_0$ 
has a surjective projection onto $2$-dimensional {\em momentum cone} 
$\check{\sigma}$ where $\sigma= {\mathcal C}(E_{h,m}) = 
\R_{\geq 0} v_3 +  \R_{\geq 0} v_4$.  The fibers of this projection 
are $1$-dimensional {\em string polytopes}.   
Since $p+q \neq 1$,  the affine toric variety 
$T_{h,m}$ does not admit a quasihomogeneous $SL(2)$-action (see also 
a remark in \cite[Section 8]{Ga08}).

\begin{remark}  
{\rm It is not easy to describe the behavior of toric degenerations 
under equivariant morphisms of spherical varieties. The simplest 
example in \ref{tor-deg} shows that toric degenerations do not 
preserve equivariant open embeddings: toric geneneration  ${\mathcal U}_0$ 
of the open orbit ${\mathcal U} \subset E_{h,m}$ is not an 
open subset in $T_{h,m}$, the corresponding birational 
morphism  ${\mathcal U}_0 \to T_{h,m}$ contracts a divisor in 
${\mathcal U}_0$.   

We remark that if $m =1$ then 
$T_{h,m}^+$ locally isomorphic to product $\A^2/\mu_p \times \A^1$. 
Therefore, toric degeneration $T_{h,m}^+$ of $E_{h,m}^+$ 
has the same type of toroidal singularity 
along the curve $C_T^+ \subset T_{h,m}^+$ as $C^+ \subset E_{h,m}^+$
However, the same is not true for the  toric degeneration $T_{h,m}^-$ 
of $E_{h,m}^-$. For instance,  if $m =1$ then  
$T_{h,1}^-$ has only a single isolated singularity, but  
singular locus  of $E_{h,1}^-$ is the whole curve 
$C^- \subset  E_{h,1}^-$. 
} 
\end{remark}

\end{document}